\documentclass[11pt]{article}

\date{April 16, 2020}

\usepackage{graphicx,authblk}
\usepackage[authoryear]{natbib}
\usepackage{fullpage}
\usepackage{amsmath}
\usepackage{amssymb}
\usepackage{amsthm}
\usepackage{graphics}
\usepackage{epsfig}
\usepackage{algorithm}
\usepackage{algorithmic}
\usepackage{mathrsfs}
\usepackage{url}
\usepackage{enumitem}
\usepackage{subfig}
\usepackage{color}
\usepackage{graphicx}
\usepackage{epstopdf}
\usepackage{tikz-qtree,tikz}
\usepackage{hyperref}
\usepackage{multicol}
\usepackage[funcfont=italic,full]{./complexity/complexity}

\newcommand{\Home}{./}

\newcommand{\ImageDir}{\Home/images}

\graphicspath{{\ImageDir/}}
\makeatletter
\def\input@path{{\ImageDir/}}
\makeatother

\usepackage[normalem]{ulem}
\usepackage{todonotes}
\usepackage{boxedminipage}

\usepackage{setspace}

\newtheorem{theorem}{Theorem}[section]

\theoremstyle{definition}
\newtheorem{definition}{Definition}
\newtheorem{assumption}{Assumption}

\newcommand{\Z}{\mathbb{Z}}
\renewcommand{\P}{\mathcal{P}}
\newcommand{\F}{\mathcal{F}}
\newcommand{\Q}{\mathbb{Q}}
\renewcommand{\E}{\mathbb{E}}
\newcommand{\J}{L}
\newcommand{\proj}{\operatorname{proj}}
\newcommand{\argmin}{\operatorname{argmin}}

\newcommand{\conv}{\operatorname{conv}}
\newcommand{\m}[1]{$#1$}

\newcommand{\T}{\mathcal{T}}

\newcommand{\Scal}{\mathcal{S}}

\newcommand{\phiLP}{\phi_{\rm LP}}
\newcommand{\uphiLP}{\underline{\phi}_{\rm LP}}
\newcommand{\ra}{\rightarrow}
\newcommand{\epi}{\operatorname{epi}}
\newcommand{\yo}{y^\omega}

\newcommand{\midd}{\;\middle|\;}
\newcommand{\RR}{\mathcal{R}}

\definecolor{Red}{RGB}{255,0,0}
\definecolor{LightCoral}{RGB}{240,128,128}
\definecolor{DeepPink}{RGB}{255,20,148}
\definecolor{LightSalmon}{RGB}{255,160,122}
\definecolor{HotPink}{RGB}{240,120,160}
\definecolor{Indigo}{RGB}{71,60,139}
\definecolor{DodgerBlue}{RGB}{30,144,255}
\definecolor{Navy}{RGB}{0,0,128}
\definecolor{DeepSkyBlue}{RGB}{0,191,255}
\definecolor{SteelBlue}{RGB}{70,130,180}
\definecolor{MediumSeaGreen}{RGB}{60,179,113}
\definecolor{DarkGreen}{RGB}{0,100,0}
\definecolor{Khaki}{RGB}{240,230,160}
\definecolor{DarkGoldenrod}{RGB}{184,134,11}

\DeclareMathOperator{\st}  {s.\!t.}
\renewcommand{\Re}{\mathbb{R}}
\renewcommand{\S}{\mathcal{S}}

\newcommand{\noprint}[1]{}

\newcounter{examples}
\newenvironment{example}[1][Example \theexamples]
{\refstepcounter{examples} 
  \begin{trivlist} \item[\hskip \labelsep {\bfseries #1}]}
  {\qed \end{trivlist}}

\makeatletter
\newenvironment{mysubequations}[1]
 {%
  \addtocounter{equation}{-1}%
  \begin{subequations}
  \renewcommand{\theparentequation}{#1}%
  \def\@currentlabel{#1}%
 }
 {%
  \end{subequations}\ignorespacesafterend
 }
\makeatother

\begin{document}

\title{A Unified Framework for Multistage and Multilevel Mixed Integer Linear
  Optimization}
 
\author[1]{Suresh Bolusani\thanks{\texttt{bsuresh@lehigh.edu}}}
\author[2]{Stefano Coniglio\thanks{\texttt{s.coniglio@soton.ac.uk}}}
\author[1]{Ted K. Ralphs\thanks{\texttt{ted@lehigh.edu}}}
\author[3]{Sahar Tahernejad\thanks{\texttt{sahar@lindo.com}}}
\affil[1]{Department of Industrial and Systems Engineering, Lehigh University, Bethlehem, PA}
\affil[2]{Department of Mathematical Sciences, University of Southampton,
  Southampton, UK}
\affil[3]{Lindo Systems, Inc., Chicago, IL 60642}
  
\maketitle

\begin{abstract}
We introduce a unified framework for the study of multilevel mixed integer linear optimization problems and multistage stochastic mixed integer linear optimization problems with recourse.
The framework highlights the common mathematical structure of the two problems and allows for the development of a common algorithmic framework.
Focusing on the two-stage case, we investigate, in particular, the nature of the value function of the second-stage problem, highlighting its connection to dual functions and the theory of duality for mixed integer linear optimization problems, and summarize different reformulations.
We then present two main solution techniques, one based on a Benders-like decomposition to approximate either the risk function or the value function, and the other one based on cutting plane generation.
\end{abstract}

\section{Introduction}\label{sec:intro}

This article introduces a unified framework for the study of \emph{multilevel
mixed integer linear optimization problems} and \emph{multistage stochastic
mixed integer linear optimization problems with recourse}. 
This unified framework
provides insights into the nature of these two
well-known classes of optimization problems, highlights their common
mathematical structure, and allows results from the wider literature devoted
to both classes to be exploited for the development of a common
algorithmic framework.

\subsection{Motivation}

Historically, research in mathematical optimization, which is arguably the
most widely applied theoretical and methodological framework for solving
optimization problems, has been primarily focused on ``idealized'' models
aimed at informing the decision process of a single
\emph{decision maker} (DM)
facing the problem of making a single set of decisions at a single point in
time
under perfect information.
Techniques for this idealized case are now well developed, with
efficient implementations widely available in off-the-shelf software.

In contrast, most real-world applications involve multiple DMs,
and decisions must be made at multiple points in time under uncertainty.
To allow for this additional complexity, a number of more sophisticated
modeling frameworks have been developed, including multistage and multilevel
optimization.
In line with the recent optimization literature, we use the term
\emph{multistage optimization} to denote the decision process of a single DM
over multiple time periods with an objective that factors in the (expected)
impact at future stages of the decisions taken at the current stage.
With the term \emph{multilevel optimization}, on the other hand,
we refer to game-theoretic
decision processes in which multiple DMs with selfish objectives make
decisions in turn, competing to optimize their own individual outcomes in the
context of settings such as, e.g., economic markets.

Because the distinction between multistage and multilevel optimization
problems appears substantial from a modeling perspective, their development has
been undertaken independently by different research communities.
Indeed, multistage problems have arisen out of the necessity to account for
stochasticity, which is done by explicitly including multiple
decision stages in between each of which the value of a random
variable is realized.
Knowledge of the precise values of the quantities that were unknown at
earlier stages allows for
so-called \emph{recourse decisions} to be made in later stages in order to
correct possible mis-steps
taken due to the lack of precise information.
On the other hand, multilevel optimization has been developed primarily to model
multi-round games (technically known as \emph{finite, extensive form games with perfect information} in
the general case and \emph{Stackelberg games} in the case of two rounds) in which the decision (or strategy) of a given player at a given
round must take into account the reactions of other players in future rounds.

Despite these distinctions, these classes of problems share an important
common structure from a mathematical and methodological perspective that makes
considering them in a single, unifying framework attractive.
It is not difficult to understand the source of this common structure---from
the standpoint of an individual DM, the complexity of the decision process
comes from uncertainty about the future.
From an analytical perspective, the methods we use for dealing with the
uncertainty arising from a lack of knowledge of the precise values of input
parameters to later-stage decision problems can also be used to address
the uncertainty arising from a lack of knowledge of the future actions of
another self-interested player.
In fact, one way of viewing the outcome of a random variable is as a
``decision'' made by a DM about whose objective function nothing is known.
Both cases require consideration of a set of outcomes arising from
either the different ways in which the uncertain future could unfold
or from the different possible actions the other players could take.
Algorithms for solving these difficult optimization problems must, either
explicitly or implicitly, rely on efficient methods for exploring this outcome
space. This commonality turns out to be more than a philosophical abstraction.
The mathematical connections between multistage and multilevel optimization
problems run deep and existing algorithms for the two cases already exhibit
common features, as we illustrate in what follows. 

\subsection{Focus of the Paper}

In the rest of
this paper, we 
address the broad class of optimization problems that we refer to from here on
by the collective name \emph{multistage
mixed integer linear optimization problems}. Such problems allow for multiple
decision stages, with decisions at each stage made by a different DM and with
each stage followed by the revelation of the value of one or more random
variables affecting the available actions at the subsequent stages. Each DM is
assumed to have their own objective function for the evaluation of the decisions
made at all stages following the stage at which they make their first decision,
including stages whose decision they do not control. Importantly, the
objective functions of different DMs may (but do not necessarily) conflict.
Algorithmically, the focus of such problems is usually on determining an
optimal decision at the first stage. At the point in time when later-stage
decisions must be made, the optimization problem faced by those DMs has the
same form as the one faced by early-stage DMs but, in it, the decisions made
at the earlier stages act as deterministic inputs. Note that, while we have
assumed different DMs at each stage, it is entirely possible to model scenarios
in which a single DM makes multiple decisions over time. From a mathematical
standpoint, this latter situation is equivalent to the case in which different DMs
share the same objective function and we thus do not differentiate these
situations in what follows.

Although the general framework we introduce applies more broadly, we focus
here on problems with two decision stages and two DMs, as well as stochastic
parameters whose values are realized between the two stages. We further
restrict our consideration to the case in which we have both continuous and
integer variables
but the constraints and objective functions
are linear.
We refer to these problems as \emph{two-stage mixed integer linear
  optimization problems} (2SMILPs).
Despite this restricted setting, the framework
can be extended to multiple stages and more general forms of constraints and
objective functions in a conceptually straightforward way (see, e.g.,
Section~\ref{sec:alternative-models}).

\section{Related and Previous Work}

In this section, we give a brief overview of related works. The literature on
these topics is vast and the below overview is not intended to be exhaustive
by any means, but only to give a general sense of work that has been done to
date. The interested reader should consult other articles in this volume for
additional background and relevant citations. 

\subsection{Applications}

Multilevel and multistage structures, whose two level/stage versions are known
as \emph{bilevel optimization problems} and \emph{two-stage stochastic
  optimization problems with recourse} (2SPRs), arise naturally in a vast array
of applications of which we touch on only a small sample here. 

In the modeling of large organizations with hierarchical decision-making
processes, such as corporations and governments, \cite{bard83} discusses a
bilevel corporate structure in which top management is the leader and
subordinate divisions, which may have their own conflicting objectives, are
the followers.
Similarly, government policy-making can be viewed from a bilevel optimization
perspective:~\cite{bard-plummer-sourie00} models a government encouraging
biofuel production through subsidies to the petro-chemical industry,
while~\cite{amouzegar-moshirvaziri99} models a central authority that sets
prices and taxes for hazardous waste disposal while polluting firms respond
by making location-allocation and recycling decisions.

A large body of work exists on \emph{interdiction problems}, which model
competitive games where two players have diametrically opposed goals and the
first player has the ability to prevent one or more of the second player's
possible activities (variables) from being engaged in at a non-zero
level. 
Most of the existing literature on
these problems has focused on variations of the well-studied \emph{network
  interdiction problem}~\cite{wollmer64,mcmasters70, ghare71, wood93,
  cormican98, israeli02, held05, janjarassuk06}, in which the lower-level DM 
is an
entity operating a network of some sort and the upper-level DM (or
interdictor) attempts to {\em interdict} 
the network as much as possible via the removal (complete or otherwise) of
portions (subsets of arcs or nodes) of the network. 
A more general case which does not involve networks (the so-called
\emph{linear system interdiction problem}) was studied in~\cite{israeli99} and
later in~\cite{DeNegre2011}.
A related set of problems involves an attacker disrupting warehouses or other
facilities to maximize the resulting transportation costs faced by the firm
(the follower)~\cite{ChScMi04,ScaparraChurch08COR,ZhaSnyRalXue16}. 
A \emph{trilevel} version of this problem involves the firm first fortifying
the facilities, then the attacker interdicting them, and finally the firm
re-allocating customers~\cite{ChSc06}. 
More abstract graph-theoretical interdiction problems in which the vertices of
a graph are removed in order to reduce the graphs' stability/clique number are
studied
in~\cite{rutenburg1994propositional,furinimaximum,coniglio2017separation}. 

Multilevel problems arise in a wide range of industries.
For instance, in the context of the electricity industry, \cite{hobbs-nelson92}
applies bilevel optimization to demand-side management while
\cite{GRIMM2016493} formulates a trilevel problem with power network expansion
investments in the first level, market clearing in the second, and redispatch
in the third. \cite{cote-marcotte-savard03,coniglioAirline} addresses the
capacity allocation and pricing problem for the airline industry. \cite{dempe-kalashnikov-mercado05} presents a model for the natural-gas
shipping industry.
A large amount of work has been carried out in the context of traffic planning
problems, including constructing road networks to maximize users'
benefits~\cite{benayed-etal92}, toll revenue
maximization~\cite{labbe98}, and hazardous material
transportation~\cite{kara-verter04}. 
For a general review on these problems and for one specialized to
price-setting, the reader is referred to~\cite{migdalas95,labbe2013bilevel}. 
More applications arise in chemical engineering and bioengineering in the
design and control of optimized systems. 
For example, \cite{clark-westerberg90} optimizes a chemical process by
controlling temperature and pressure (first-stage actions) 
where the system (second stage) reaches an equilibrium as it naturally
minimizes the Gibbs free energy. 
\cite{burgard-pharkya-maranas03} develops gene-deletion strategies (first
stage) to allow overproduction of a desired chemical by a cell (second stage). 
In the
area of telecommunication networks, bilevel optimization
has been used for modeling
the behavior of a networking communication protocol (second-level problem)
which the network operator, acting as first-level DM, can influence but not
directly control. 
The case of routing over a TCP/IP network is studied
in~\cite{amaldi2013single,amaldi2013network,amaldi2013energy,amaldi2014maximum}.  

The literature on game theory features many works on bilevel
optimization
problems naturally arising from the computation of \emph{Stackelberg
  equilibria} in different settings. 
Two main variants of the Stackelberg paradigm
are typically considered: one in which the followers can observe the action
that the leader draws from its commitment and, therefore, the commitment is in
pure strategies~\cite{von2010market}, and one in which the followers cannot do
that directly and, hence, the leader's commitment can be in mixed
strategies~\cite{conitzer2006computing,von2010leadership}. 
While most of the works focus on the case with a single leader and a single
follower (which leads to a proper bilevel 
optimization
problem), some work has been done on the case with more than two players:
see~\cite{conitzer2011commitment,basilico2016methods,basilico2017bilevel,coniglio2017pessimistic,basilico2020bilevel,marchesi2018leadership,castiglioni2019leadership,coniglio2020computing}
for the \emph{single}-leader \emph{multi}-follower case,
\cite{smith2014multidefender,lou2015equilibrium,laszka2016multi,lou2017multidefender,gan2018stackelberg}
for the \emph{multi}-leader \emph{single}-follower case,
or~\cite{castiglioni2019bealeader,pang2005quasi,leyffer2010solving,kulkarni2014shared} 
for the \emph{multi}-leader \emph{multi}-follower case. 
Practical applications are often found in \emph{security games}, which
correspond to competitive situations where a defender (leader) has to allocate
scarce resources to protect valuable targets from an attacker
(follower)~\cite{paruchuri2008playing,KiekintveldJTPOT09,an2011guards,tambe2011security}. 
Other practical applications are found in, among others, inspection
games~\cite{avenhaus1991inspector} and mechanism
design~\cite{sandholm2002evolutionary}.
The works on the computation of a {\em correlated equilibrium}~\cite{celli2019computing} as well those on Bayesian persuasion~\cite{celli2020bayesian}, where a leader affects the behavior of the follower(s) by a {\em signal}, also fall in this category.

Finally, there are deep connections between bilevel optimization and the
algorithmic decision framework that drives branch and bound itself,
and it is likely that the study of bilevel optimization problems may lead to
improved methods for solving single-level optimization problems. 
For example, the problem of determining the disjunction whose imposition
results in the largest bound improvement within a branch-and-bound framework
and the problem of determining
the
maximum bound-improving inequality are
themselves bilevel optimization
problems~\cite{Mahajan2009,MahRal10,coniglio2015generation}.
The same applies in $n$-ary branching when one looks for a branching decision
leading to the smallest possible number of child
nodes~\cite{lodi2009interdiction,san2019new}. 

Multistage problems and, in particular, two-stage stochastic optimization
  problems with recourse, arise in an equally wide array of application areas,
  including scheduling, forestry, pollution control, telecommunication and
  finance. 
  \cite{grass2016two} surveys literature, applications, and methods for
  solving disaster-management problems arising in the humanitarian context. 
  \cite{gupta2007lp} addresses network-design problems where, in the second
  stage, after one of a finite set of scenarios is realized, additional edges of
  the network can be bought, and provides constant-factor approximation
  algorithms. 
  A number of works address the two-stage stochastic optimization with
  recourse version of classical combinatorial optimization problems: among
  others, \cite{dhamdhere2005two} considers the spanning-tree problem,
  \cite{katriel2007commitment} the matching problem, and
  \cite{gendreau1996stochastic,gortz2012stochastic} the vehicle routing
  problem. 
  For references to other areas of applicability, see the 
  books~\cite{birge2011introduction,kall2010stochastic} and, in
  particular,~\cite{wallace2005applications}. 

\subsection{Algorithms}

The first recognizable formulations for bilevel optimization problems were
introduced in the 1970s in~\cite{bracken-mcgill73} and this is when the term
was also first coined.
Beginning in the early 1980s, these problems attracted increased interest.
\cite{vicente-calamai94} provide a large bibliography of the early developments.

There is now a burgeoning literature on continuous bilevel linear
optimization, but it is only in the past decade that work on the discrete case
has been undertaken in earnest by multiple research groups.
\cite{moore90} was the first to introduce a framework for general integer
bilevel linear optimization and to suggest a simple branch-and-bound
algorithm. The same authors also proposed a more specialized algorithm for
binary bilevel optimization problems in~\cite{bard92}.
Following these early works, the focus shifted primarily to various special
cases, especially those in which the lower-level problem has the integrality
property.
\cite{dempe01} considers a special case characterized by continuous
upper-level variables and integer lower-level variables and uses a cutting
plane approach to approximate the lower-level feasible region (a somewhat
similar approach is adopted in~\cite{Dempe2017} for solving a bilinear mixed
integer bilevel problem with integer second-level variables). 
\cite{wen90} consider the opposite case, where the lower-level problem is a
linear optimization problem and the upper-level problem is an integer
optimization problem, using  linear optimization duality to derive exact and
heuristic solutions. 

The publication of a general algorithm for pure integer problems
in~\cite{DeNRal09} (based on the groundwork laid in a later-published
dissertation~\cite{DeNegre2011}) spurred renewed interest in developing
general-purpose algorithms. The evolution of work is summarized in
Table~\ref{table:previousWork}, which indicates the types of variables
supported in both the first and second stages (C indicates continuous, B
indicates binary, and G indicates general integer). The aforementioned network
interdiction problem is a special case that continues to receive significant
attention, since tractable duality conditions exist for the lower-level
problem~\cite{wollmer64, mcmasters70, ghare71, wood93, cormican98, israeli02,
  held05, janjarassuk06}.

\begin{table}[]
\begin{center}
\begin{tabular}{|l|l|l|}
\hline
\multicolumn{1}{|l|}{Citation}
& \multicolumn{1}{l|}{Stage 1 Variable Types} & \multicolumn{1}{l|}{Stage 2
  Variable Types} \\ \hline
\cite{wen90} & B & C \\
\cite{bard92} & B & B\\
\cite{faisca-etal07} & B, C & B, C \\
\cite{saharidis-ierapetritou08} & B, C & B, C \\
\cite{garcesetal09} & B & C \\
\cite{DeNRal09}, \cite{DeNegre2011} & G & G \\
\cite{koppe10} & G or C & G \\
\cite{baringo-conejo12} & B, C & C \\
\cite{xuwang14} & G & G, C \\
\cite{zengan14} & G, C & G, C \\
\cite{caramiamari15} & G & G \\
\cite{capraraetal16} & B & B \\
\cite{HemSmi16} & B, C & B, C \\
\cite{TahRalDeN16} & G, C & G, C \\
\cite{wangxu17} & G & G \\
\cite{lozanosmith17} & G & G, C \\
\cite{fischettietal17a}, \cite{fischettietal17b} & G, C & G, C \\
\hline
\end{tabular}
\end{center}
\caption[Previous Work]{Evolution of algorithms for bilevel
  optimization}\label{table:previousWork} 
\end{table}

As for the case of multistage stochastic optimization, the two-stage linear 
stochastic optimization problem with recourse in which both the first- and second-stage problems
contain only continuous variables has been well studied both theoretically and
methodologically.
\cite{birge2011introduction,kall2010stochastic} survey the related literature.
The integer version of the problem was first considered in the early 1990s by
\cite{louveaux1993stochastic} for the case of two-stage problems with
\emph{simple integer recourse}. 
Combining the methods developed for the linear version with the
branch-and-bound procedure, \cite{laporte1993integer} proposed an algorithm
known as the \emph{integer L-shaped method} where the first-stage problem
contains only binary variables. 
Due to
the appealing
structural properties of a (mixed) binary integer optimization problem, a
substantial amount of literature since then has been considering the case of a
two-stage stochastic problem with (mixed) binary variables in one or both
stages~\cite{sen2005c,sherali2006solving,gade2012decomposition}. 
The case of two-stage problems with a pure integer recourse has also been
frequently visited, see~\cite{schultz1998solving,kong2006two}. 
It must be noted that methods such as these, which are typically developed for
special cases, often rely on the special structure of the 
second-stage problem, thus being often not applicable to the two-stage problem
with mixed integer restrictions. 
Algorithms for stochastic optimization problems with integer recourse were
proposed by~\cite{caroe1998shaped} and~\cite{sherali2002modification}. 

\section{Setup and preliminaries \label{sec:setup}}

The defining feature of a multistage optimization problem is that the values
of the \emph{first-stage} variables (sometimes called \emph{upper-level} or
\emph{leader} variables in the bilevel optimization literature) must be
(conceptually) fixed without explicit knowledge of future events, the course
of which can be influenced by the first-stage decision itself. Due to this
influence, the perceived ``value'' of the first-stage decision must
take into account the effect of this decision on the likelihood of occurrence
of these future events.

More concretely, the first-stage DM's overall objective is to minimize the sum
of two terms, the first term representing the immediate cost of implementation
of the first-stage solution and the second term representing the desirability
of the first-stage decision in terms of its impact on the decisions taken at
later stages.
The general form of a two-stage mixed integer linear optimization
  problem is then
\begin{equation} \label{eqn:2SMILP} \tag{2SMILP}
\min_{x \in \P_1 \cap X} \left\{ c x + \Xi(x)\right\},
\end{equation}
where
\begin{equation*}
\P_1 = \left\{x\in \Re^{n_1} \mid A^1 x \geq b^1 \right\}
\end{equation*}
is the \emph{first-stage feasible region}, with $A^1 \in \mathbb{Q}^{m_1\times
  n_1}$ and $b^1 \in \mathbb{Q}^{m_1}$ defining the associated linear
constraints, and $X = \Z^{r_1}_+ \times \mathbb{Q}^{n_1-r_1}_+$ representing
integrality, rationality, and non-negativity requirements on the first-stage variables, denoted
by $x$. Note that we require the continuous variables to take on rational
values in order to ensure that the second-stage problem defined
in~\eqref{eqn:second-stage-problem} below has an optimal value that is attainable
when that value is finite. In practice, solvers always return such solutions,
so this is a purely technical detail.
The linear function $cx$ with $c \in \mathbb{Q}^{n_1}$ is the aforementioned
term that reflects the immediate cost of implementing the first-stage solution.
The function $\Xi: \Q^{n_1} \rightarrow \Q\cup\{\pm \infty\}$ is the
\emph{risk function}, which takes only rational input for technical reasons
discussed later.
$\Xi$ is the aforementioned term representing the first-stage DM's evaluation
of the impact of a given choice for the value of the first-stage variables on
future decisions.
Similar concepts of risk functions have been employed in many different
application domains and will be briefly discussed in
Section~\ref{sec:alternative-models}.
To enable the development of a practical methodology for the solution of these
problems, however, we now define the specific class of functions we consider.

\subsection{Canonical Risk Function}

Our canonical risk function is a generalization of the risk function
traditionally used in defining 2SPRs. As usual, let us now introduce a random
variable $U$ over an outcome space $\Omega$ representing the set of possible
future scenarios that could be realized between the making of the first- and
second-stage decisions. The values of this random variable will be input
parameters to the so-called \emph{second-stage problem} to be defined below.

As is common in the literature on 2SPRs, we assume that $U$ is discrete, i.e.,
that the outcome space $\Omega$ is finite, so that $\omega \in \Omega$
represents which of a finite number of explicitly enumerated scenarios is
actually realized. In practice, this assumption is not very restrictive, as
one can exploit any algorithm for the case in which $\Omega$ is assumed finite
to solve cases where $\Omega$ is not (necessarily) finite by utilizing a
technique for discretization, such as \emph{sample average approximation}
(SAA)~\cite{shapiro2003monte}. As $U$ is discrete in this work, we can
associate with it a probability distribution defined by $p \in \Re^{|\Omega|}$
such that $0 \leq p_\omega \leq 1$ and $\sum_{\omega \in \Omega} p_\omega = 1$.

With this setup, the canonical risk function for $x \in \Q^{n_1}$ is
\begin{equation} \label{eqn:risk-function} \tag{RF}
\Xi(x) = \E \left[\Xi_\omega(x)\right] = \sum_{\omega \in \Omega} p_\omega
\Xi_\omega(x), 
\end{equation}
where $\Xi_\omega(x)$ is the \emph{scenario risk function}, defined as
\begin{equation} \label{eqn:second-stage-function} \tag{2SRF}
\Xi_\omega(x) =  \min \left\{ d^1 y^{\omega} \midd y^{\omega} \in\argmin
\{d^2y \mid y \in \P_2(b^2_\omega - A^2_\omega x) \cap Y\}
\right\};
\end{equation}
the set
\begin{align*}
\P_2(\beta) &= \left\{y\in \Re^{n_2}\mid G^2y \geq \beta \right\}
\end{align*}
is one member of a family of polyhedra that is parametric w.r.t. the
right-hand side vector $\beta \in \Re^{m_2}$ and represents the second-stage
feasibility conditions; and $Y = \Z^{r_2}_+ \times \mathbb{Q}^{n_2-r_2}_+$
represents the second-stage integrality and non-negativity requirements.
The deterministic input data defining $\Xi_\omega$ are $d^1,d^2 \in
\mathbb{Q}^{n_2}$ and $G^2\in\Q^{m_2\times n_2}$.
$A^2_\omega \in\Q^{m_2\times n_1}$ and $b^2_\omega \in \mathbb{Q}^{m_2}$
represent the realized values of the random input parameters in scenario
$\omega \in \Omega$, i.e., $U(\omega) = (A^2_\omega, b^2_\omega)$.

As indicated in~\eqref{eqn:risk-function}
and~\eqref{eqn:second-stage-function}, the inner optimization problem faced by
the second-stage DM is parametric only in its right-hand side,
which is determined jointly by the value $\omega$ of the random variable $U$
and by the chosen first-stage solution. It will be useful in what follows to
define a family of \emph{second-stage} optimization problems
\begin{equation} \tag{SS} \label{eqn:second-stage-problem}
  \inf \left\{ d^2 y \midd y \in \P_2(\beta) \cap Y \right\}
\end{equation}
that are parametric in the right-hand side $\beta \in \Re^{m_2}$ (we use
``$\inf$'' instead of ``$\min$'' here because, for $\beta \not\in
\mathbb{Q}^{m_2}$, the minimum may not exist). By further defining
\begin{align*}
  \beta^\omega(x) = b^2_\omega - A^2_\omega x
\end{align*}
to be the parametric right-hand side that arises when the chosen first-stage
solution is $x \in X$ and the realized scenario is $\omega \in \Omega$, we can
identify the member of the parametric family defined
in~\eqref{eqn:second-stage-problem} in scenario $\omega \in \Omega$ when the
chosen first-stage solution is $x \in X$ as that with feasible region
$\P_2(\beta^\omega(x)) \cap Y$.
Associated with each $x \in \Q^{n_1}$ and $\omega \in
\Omega$ is the set of all alternative optimal solutions to the second-stage
problem~\eqref{eqn:second-stage-problem} (we allow for $x \not\in X$ here
because such solutions arise when solving certain relaxations), called the
\emph{rational reaction set} and denoted by
\begin{equation*}
\RR^\omega(x) = \argmin\{d^2y \mid y \in \P_2(b^2_\omega - A^2_\omega x) \cap
Y\}. 
\end{equation*}
For a given $x \in \Q^{n_1}$, $\RR^\omega(x)$ may be empty if
$\P_2(b^2_\omega - A^2_\omega x)\cap Y$ is itself empty or if the second-stage
problem is unbounded (we assume in Section~\ref{sec:assumptions} that this
cannot happen, however).

When $|\RR^\omega(x)| > 1$, the second-stage DM can, in principle, choose which
alternative optimal solution to implement. We must therefore specify in the
definition of the risk function a rule by which to choose one of the
alternatives. According to our canonical risk
function~\eqref{eqn:risk-function} and the corresponding scenario risk
function~\eqref{eqn:second-stage-function}, the rule is to choose, for each
scenario $\omega \in \Omega$, the alternative optimal solution that minimizes
$d^1 y^\omega$, which corresponds to choosing the collection
$\{y^\omega\}_{\omega \in \Omega}$ of solutions to the individual scenario
subproblems that minimizes
\begin{align*}
  d^1 \left(\sum_{\omega \in \Omega} p_\omega y^\omega \right).
\end{align*}
This is known as the \emph{optimistic} or \emph{semi-cooperative} case in the
bilevel optimization literature, since it corresponds to choosing the
alternative that is most beneficial to the first-stage DM.
Throughout the paper, we consider this case unless otherwise specified.
In Section~\ref{sec:alternative-models}, we discuss other forms of risk
function.  

Because of the subtleties introduced above, there are a number of ways one
could define the ``feasible region'' of~\eqref{eqn:2SMILP}.
We define the \emph{feasible region} for scenario~$\omega$ (with respect to both
first- and second-stage variables) as 
\begin{equation*} \label{eqn:F-omega} \tag{FR}
\F^\omega = \{(x,y^\omega) \in X \times Y \mid x \in \P_1, y^\omega \in
\RR^\omega(x)\} 
\end{equation*}
and members of $\F^\omega$ as \emph{feasible solutions} for scenario $\omega$.
Note that this definition of feasibility
does not prevent having
$(x, y^\omega) \in \F^\omega$ but $d^1 y^\omega > \Xi_\omega(x)$. This
will not cause any serious difficulties, but is something to keep in mind.  

We can similarly define the feasible region with respect to just the
first-stage variables as 
\begin{equation}\label{eqn:f1}\tag{FR-Proj}
\F^1 = \bigcap_{\omega \in \Omega} \proj_x(\F^\omega).
\end{equation}
Since $\Xi(x) = \infty$ for $x \in \Q^{n_1}$ if the feasible region
$\P_2(\beta^\omega(x)) \cap Y$ of the second-stage
problem~\eqref{eqn:second-stage-problem} is empty for some $\omega \in
\Omega$,
we have that, for $x \in \P_1 \cap X$, the following are
  equivalent:
\begin{equation*}
x \in \F^1 \Leftrightarrow x \in \bigcap_{\omega \in \Omega}
\proj_x(\F^\omega) \Leftrightarrow \RR^\omega(x) \not= \emptyset\; \forall
\omega \in \Omega  \Leftrightarrow \Xi(x) < \infty. 
\end{equation*}

Finally, it will be convenient to define $\P^\omega$ to be the feasible region
of the relaxation of the deterministic two-stage problem under scenario
$\omega \in \Omega$ that is obtained by dropping the optimality requirement for the
second-stage variables $y^\omega$, as well as any integrality restrictions.
Formally, we have:
$$
\P^\omega = \left\{(x,y^{\omega}) \in \Re^{n_1 + n_2}_+ \mid x \in \P_1,
y^{\omega} \in \P_2(b^2_\omega- A^2_\omega x) \right\}. 
$$
Later in Section~\ref{sec:algorithms}, we will use these sets to define a
relaxation for the entire problem that will be used as the basis for the
development of a branch-and-cut algorithm. 

\subsection{Technical Assumptions \label{sec:assumptions}}

We now note the following assumptions made in the remainder of the paper. 
\begin{assumption} \label{as:boundedness}
$\P^\omega$ is bounded for all $\omega \in \Omega$.
\end{assumption}
\noindent This assumption, which is made for ease of presentation and can be
relaxed, results in the boundedness 
of~\eqref{eqn:2SMILP}.
\begin{assumption} \label{as:setJ}
All first-stage variables with at least one non-zero coefficient in
the second-stage problem (the so-called \emph{linking variables}) are integer, i.e.,
\begin{equation*}
\J = \left\{i \in \{1, \dots, n_1\}\midd a^\omega_i \neq 0 \textrm{ for
  some } \omega \in \Omega \right\}\subseteq\left\{1,...,r_1\right\},
\end{equation*}
where $a^\omega_i$ represents the $i^{\textrm{th}}$ column of matrix
$A^2_\omega$.
\end{assumption}
These two assumptions together guarantee that an optimal solution exists
whenever \eqref{eqn:2SMILP} is feasible~\cite{vicente96}.
It also guarantees that
the convex hull of $\F^\omega$ is a polyhedron, which is important for
the algorithms we discuss later.
Note that,
due to the assumption of optimism, we can assume w.l.o.g. that all first-stage
variables are linking variables by simply interpreting the non-linking
variables as belonging to the second stage.
While this may seem conceptually inconsistent with the intent of the original
model, it is not difficult to see that the resulting model is
\emph{mathematically} equivalent, since these variables do not affect the second-stage
problem and thus, the optimistic selection of values for those variables will
be the same in either case.

Before closing this section, we remark that, in this article, we do not allow
second-stage variables in the first-stage constraints. 
While this case can be handled with techniques similar to those we describe in
the paper from an algorithmic perspective, it does require a more complicated
notation which, for the sake of clarity, we prefer not to adopt.
Detailed descriptions of algorithms for this more general case in the bilevel
setting are provided in~\cite{TahRalDeN16,BolRal20}.

\subsection{Alternative Models \label{sec:alternative-models}}

\paragraph{Alternative Form of~\eqref{eqn:2SMILP}.}

For completeness, we present here an alternative form of~\eqref{eqn:2SMILP} that
is closer to the traditional form in which bilevel optimization problems are
usually specified in the literature. Adopting the traditional notation,
\eqref{eqn:2SMILP} can be alternatively written as
\begin{equation} \label{eqn:2SMILP-standard} \tag{2SMILP-Alt}
\begin{aligned} 
  & \min_{x,\{y^\omega\}_{\omega \in \Omega}} && cx + \sum_{\omega \in \Omega}
  p_\omega d^1 y^\omega\\ 
  & \st  && A^1 x  \geq b^1\\
             &&& x  \in X\\
             &&& \hspace{-.3cm}\left. \begin{array}{llllr}
                        y^\omega \in & \argmin_{y} & d^2 y \\
                        & \st & A^2_\omega x + G^2 y \geq b^2_\omega\\
                        && y \in Y  
                      \end{array} \right\} \forall \omega \in \Omega.
\end{aligned}
\end{equation}
Note that, in the first stage, the minimization is carried out with respect to
both $x$ and $y$. This again specifies the optimistic case discussed earlier,
since the above formulation requires that, for a given $x \in X$, we select
$\{y^\omega\}_{\omega \in \Omega}$ such that
\begin{align*}
  d^1 \left(\sum_{\omega \in \Omega} p_\omega y^\omega \right)
\end{align*}
is minimized. 

\paragraph{Pessimistic Risk Function.}

As already pointed out, the canonical
risk function defined in~\eqref{eqn:risk-function} assumes the
\emph{optimistic case}, since it encodes selection of the alternative optimal
solution to the second-stage problem that is most beneficial to the
first-stage DM. This is the case we focus on in the remainder of the paper.
The pessimistic case, on the other hand, is easily modeled by defining the
scenario risk function
to be
\begin{equation*} 
    \Xi_\omega(x) =  \max \left\{ d^1 y^{\omega} \midd y^{\omega} \in\argmin
    \{d^2 y \mid y \in \P_2(\beta^\omega(x) \cap Y)\} \right\}. 
  \end{equation*}
We remark that, while the optimistic and pessimistic cases may coincide in
some cases (e.g., when
\eqref{eqn:second-stage-problem} admits a single optimal solution for every
$x$),
this coincidence is rarely observed in practice and would be hard to detect in
any case. In general, the pessimistic case is more difficult to solve, though
the algorithms discussed in Section~\ref{sec:algorithms} can be modified to
handle it.

\paragraph{Recursive Risk Functions.}

Although we limit ourselves to problems with two stages in this paper, we
briefly mention that more general risk functions can be defined by recursively
defining risk functions at earlier stages in terms of later-stage risk
functions. This is akin to the recursive definition of the cost-to-go
functions that arise in stochastic dynamic programming
(see~\cite{bertsekas2017dynamic}). With such recursive definitions, it is
possible to generalize much of the methodology described here in a relatively
straightforward way, though the algorithm complexity grows exponentially with
the addition of each stage. It is doubtful exact algorithms can be made
practical in such cases.

\subsubsection{Other Risk Functions.}

Other forms of risk function have been used
in the literature, especially in finance. 
In robust optimization, for example, one might consider a risk function of the
form 
\begin{equation*} 
\Xi(x) = \max_{\omega \in \Omega} \left\{\Xi_\omega(x)\right\},
\end{equation*}
which models the impact on the first-stage DM of the worst-case second-stage
realization of the random variables. 
A popular alternative in finance applications that is slightly less
conservative is the \emph{conditional value at risk}, the expected value taken
over the worst $\alpha$-percentile of outcomes~\cite{rockafellar2000optimization,uryasev2000conditional}.
While it is possible
to incorporate such risk functions into the general algorithmic framework we
present here, for the purposes of limiting the scope of the discussion, we
focus herein only on risk functions in the canonical
form~\eqref{eqn:risk-function}.

\subsection{Related Classes}\label{subsec:related_classes}

With $\Xi$ defined as in~\eqref{eqn:risk-function}, the
problem~\eqref{eqn:2SMILP} generalizes several well-known classes of
optimization problems. 

\subsubsection{Single-Stage Problems.}
When $d^1 = d^2$ and $|\Omega| = 1$, the two stages of~\eqref{eqn:2SMILP} can
be collapsed into a single stage and the problem reduces to a traditional
mixed integer linear optimization problem (MILP). 
It is natural that algorithms for~\eqref{eqn:2SMILP} rely heavily on solving
sequences of related single-stage MILPs and we discuss parametric versions of
this class in later sections. For continuity, we utilize the notation for the
second-stage variables and input data throughout.
The case of $r_2 = 0$ (in which there are no integer variables) further
reduces to a standard linear optimization problem (LP). 

\subsubsection{Bilevel Problems.}

When $|\Omega|=1$ and assuming that we may have $d^1 \not=
d^2$, \eqref{eqn:2SMILP} takes the form of a \emph{mixed integer bilevel
linear optimization problem} (MIBLP). Dropping the scenario super/subscript
for simplicity, this problem is more traditionally written as
\begin{equation}\label{eqn:MIBLP} \tag{MIBLP}
  \min_{x \in \P_1 \cap X, y \in Y} \left\{c x + d^1 y \,\Big|\,
    y \in  \underbrace{\argmin\{d^2 y \mid y \in \P_2(b^2 - A^2 x) \cap
    Y\}}_{\RR(x)} 
  \right\}.
\end{equation}
Note that this formulation implicitly specifies the optimistic case, since if
$\RR(x)$ is not a singleton, it requires that among the alternative optima,
the solution minimizing $d^1 y$ be chosen. In this setting, the \emph{bilevel
risk function} can be written as
\begin{equation*} 
 \begin{aligned}
\Xi(x) &=  \min \left\{ d^1 y \mid y \in \RR(x) \right\}.
\end{aligned}
\end{equation*}

\subsubsection{Two-Stage Stochastic Optimization Problems with Recourse.}
When $d^1=d^2$, either the inner or the outer minimization
in~\eqref{eqn:second-stage-function} is redundant and \eqref{eqn:2SMILP} takes
the form of a \emph{two-stage stochastic mixed integer linear optimization
  problem with recourse}. 
In this case, for each scenario $\omega \in \Omega$ we can write the scenario risk function more simply as
\begin{equation*} 
\Xi_\omega(x) =  \min \big\{ d^1 y^\omega \mid y^\omega \in \P_2(b^2_\omega -
A^2_\omega x) \cap Y \big\}. 
\end{equation*}
The second-stage solution $y^\omega$ corresponding to scenario $\omega \in
\Omega$ is usually called the \emph{recourse decision}. 
These problems involve a single DM optimizing a single objective function, but
capable of controlling two sets of variables: the first-stage {\em
  here-and-now} variables $x$ and the second-stage \emph{wait-and-see}
variables $y^\omega$, whose value is set after observing the realization of
the random event $\omega$.

\subsubsection{Zero-Sum and Interdiction Problems.}

For $d^1 = -d^2$ (and typically, $|\Omega| = 1$),
\eqref{eqn:2SMILP} subsumes the case of \emph{zero-sum problems},
which model competitive games in which two players have exactly opposing goals.
An even more specially-structured subclass of zero-sum problems are
\emph{interdiction problems}, in which the first-stage variables are in
one-to-one correspondence with those of the second stage and represent the
ability of the first-stage DM to ``interdict'' (i.e., forcing to take value
zero) individual variables of the second-stage DM.
Formally, the effect of interdiction can be modeled using a variable
upper-bound constraint 
\begin{equation*}
y\leq u(e-x)
\end{equation*}
in the
second-stage problem, where $u\in \mathbb{R}^n$ is a vector of natural upper
bounds on the vector of variables $y$ and~$e$ is an $n$-dimensional column
vector of ones (here, $n = n_1 = n_2$).
Formally, the \emph{mixed integer interdiction problem}
is
\begin{equation*} 
\max_{x\in \P_1 \cap X} \min_{y\in \P_2(x) \cap Y} d^2y 
\end{equation*}
where (abusing notation slightly), we have
\begin{align*}
\P_2(x) & = \left\{y \in \Re^{n_2} \mid G^2 y \geq b^2, y \leq u(e-x)\right\}.
\end{align*}

\section{Computational Complexity}

Within the discrete optimization community, the framework typically used for
assessing problem complexity is based primarily on the well-known theory of
$\NPcomplexity$-completeness, which has evolved from the foundational work
of~\cite{cook1971complexity}, \cite{karp1975computational},
and~\cite{garey79}.
This has lead to the ubiquitous practice of classifying optimization problems as
being either in the class $\Pcomplexity$ or the class $\NPcomplexity$-hard,
the latter being an all-encompassing and amorphous class that includes
essentially all optimization problems not known to be polynomially solvable.
This categorization lacks the refinement necessary for consideration of classes
such as those described in this article. 
It is indeed easy to show that multistage optimization problems are
$\NPcomplexity$-hard in general~\cite{calamai94, Jeroslow1985, benayed90,
  hansen92}, but this merely tells us that these problems are not in
$\Pcomplexity$ (assuming $\Pcomplexity \not= \NPcomplexity$), which is not
surprising. 
What we would really like to know is for which complexity class (the decision
versions of) these problems are \emph{complete}.

In the presence of a hierarchical structure with $k$ levels (and when $\Omega$
is a singleton), the natural complexity class to consider is
$\Sigma_k^\Pcomplexity$, i.e., the $k^{\textrm{th}}$ level of the
\emph{polynomial hierarchy}.
From an optimization perspective, this hierarchy (originally introduced
in~\cite{stockmeyer77}) is a scheme for classifying multilevel decision
problems beyond the usual classes $\Pcomplexity$ and $\NPcomplexity$. 
The class $\Pcomplexity$ (which contains all decision problems that can be
solved in polynomial time) occupies the $0^{\textrm{th}}$ level,
also known as $\Sigma_0^\Pcomplexity$.
The first level, $\Sigma_1^\Pcomplexity$, is the class also known as
$\NPcomplexity$, which consists of all problems for which there exists a
certificate verifiable in polynomial time or, equivalently, all problems that
can be solved in non-deterministic polynomial time.
The $k^{\textrm{th}}$ level, $\Sigma_k^\Pcomplexity$, contains all problems with
certificates that can be be verified in polynomial time (equivalently, all
problems solvable in non-deterministic polynomial time), assuming the existence of an oracle
for solving problems in the class $\Sigma_{k-1}^\Pcomplexity$.
While it is clear that $\Sigma_{k}^\Pcomplexity \subseteq
\Sigma_{\ell}^\Pcomplexity$ for any $k, \ell \in \mathbb{N} \cup \{0\}$ with
$k \leq \ell$, $\Sigma_{k}^\Pcomplexity \subset \Sigma_{\ell}^\Pcomplexity$ is
conjectured to hold for all $k, \ell \in \mathbb{N} \cup \{0\}$ with $k <
\ell$ (the well-known $\Pcomplexity \neq \NPcomplexity$ conjecture is a
special case).
It is also known that $\Sigma_{k}^\Pcomplexity = \Sigma_{k+1}^\Pcomplexity$ would
imply $\Sigma_{k}^\Pcomplexity = \Sigma_{\ell}^\Pcomplexity$ for all $\ell \geq
k+1$, which would cause the polynomial hierarchy to collapse to level $k$ (for
$k = 0$, we would have $\Pcomplexity = \NPcomplexity$). 
The notions of completeness and hardness commonly used for $\NPcomplexity$
translate directly to $\Sigma_k^\Pcomplexity$. 
A proof that $k$-level optimization problems with binary variables, linear
constraints, and linear objective functions are hard for
$\Sigma_k^\Pcomplexity$ is contained in~\cite{Jeroslow1985}. 
Such result suffices to show that the multistage problems with $k$ stages
treated in this paper are (in their optimization version)
$\Sigma_k^\Pcomplexity$-hard (and those with $k=2$ stages are
$\Sigma_2^\Pcomplexity$-hard). 
A compendium of $\Sigma_2^\Pcomplexity$-complete/hard problems, somewhat
similar in spirit to~\cite{garey79}, can be found
in~\cite{schaefer2002completeness}, with more recent updates available online. 

For the case of two-stage stochastic
optimization problems with recourse with linear constraints, linear objective
functions, and mixed integer variables, the assumption of a finite
outcome space $\Omega$ of either fixed or polynomially bounded size suffices
to guarantee that the decision version of such a problem is $\NPcomplexity$-complete. 
Indeed, when $|\Omega|$ is considered a constant or is bounded by a polynomial
in the total number of variables and constraints, one can directly introduce
a block-structured reformulation of the problem with one block per scenario
$\omega \in \Omega$ that contains the coefficients of the constraints that
$y^\omega$ should satisfy (we discuss such a reformulation in
Section~\ref{sec:reformulations}).
As such reformulation is of polynomial size, solutions to the corresponding
optimization problem can clearly be certified in polynomial time by checking
that they satisfy all the polynomially-many constraints featured in the
formulation, which, in turn, implies that the problem belongs to
$\NPcomplexity$.
When 
the outcome space $\Omega$ is continuous, the problem becomes
$\#\Pcomplexity$-hard in general~\cite{Dyer2006,Hanasusanto2016}. 
While a single sample average approximation problem with a finite or
polynomially-bounded number of samples can be used to approximate a continuous
problem by solving a single discrete optimization problem of polynomial size,
\cite{Hanasusanto2016} shows that even finding an approximate solution using
the SAA method is $\#\Pcomplexity$-hard.
New results on the complexity of 2SPRs featuring a double-exponential algorithm can be found in~\cite{klein2019complexity}.

\section{Duality and the Value Function \label{sec:duality}}

Virtually all algorithms for the exact solution of optimization problems produce a
proof of optimality that depends on the construction of a solution to a {\em
  strong dual}.
Although the duality theory for MILPs is not widely known, the most effective
algorithms for solving MILPs (which are variants of the well-known
branch-and-bound algorithm) \emph{do} produce a solution to a certain dual
problem. 
A natural approach to solving
\eqref{eqn:2SMILP} is therefore to embed the production of the ``dual proof''
of optimality of the second-stage problem~\eqref{eqn:second-stage-problem}
into the formulation of the first-stage problem, reducing the original
two-stage problem to a traditional single-stage optimization problem.

The reformulations and algorithmic approaches that we present in
Sections~\ref{sec:reformulations} and~\ref{sec:algorithms} all use some
variant of this strategy.
In particular, the algorithms we describe are based on iteratively
strengthening an initial relaxation in a fashion reminiscent of many
iterative optimization algorithms. 
The strengthening operation essentially consists of the dynamic construction
of both a proof of optimality of the second-stage problem and of corresponding
first- and second-stage solutions. 

In the remainder of the section, we introduce the central concepts of a
duality theory for mixed integer linear optimization problems (and more
general discrete optimization problems), emphasizing its connection to
solution methods for \eqref{eqn:2SMILP}.
This introduction is necessarily brief and we refer the reader
to~\cite{GuzRal07,HasRal14-1,HasRal14} for more details specific to the
treatment here and to~\cite{wolsey81,williams96-2} for earlier foundational work
on IP duality. Although the ``dual problem'' is usually a fixed
(non-parametric) optimization problem associated with a fixed (non-parametric)
``primal problem,'' the typical concepts of duality employed in constructing
dual proofs of optimality and in designing solution algorithms inherently
involve parametric families of optimization problems. This makes the tools
offered by this theory particularly suitable for employment in this setting. To
preserve the connection with the material already introduced, we consider the
family of MILPs parameterized on the right-hand side $\beta \in \Re^{m_2}$
that was introduced earlier as~\eqref{eqn:second-stage-problem} and use the
same notation. We reproduce
it here for convenience:
\begin{equation} \tag{SS} 
  \inf \left\{ d^2 y \midd y \in \P_2(\beta) \cap Y \right\},
\end{equation}
where
\begin{align*}
\P_2(\beta) &= \left\{y\in \Re^{n_2}\mid G^2y \geq \beta \right\}
\end{align*}
and $\beta \in \Re^{m_2}$ is the input parameter. When we want to refer to a
(fixed) generic instance in this parametric family, the notation $b$ will be
used to indicate a fixed (but arbitrary) right-hand side. We also refer to
specific right-hand sides arising in the solution of~\eqref{eqn:2SMILP} using
the notation defined earlier.

\subsection{Value Functions}

Among possible notions of duality, the one most relevant to the
development of optimization algorithms is one that also has an intuitive
interpretation in terms of familiar economic concepts. This theory rests
fundamentally on an understanding of the so-called \emph{value function},
which we introduce below. The value function of an MILP has been studied by a
number of authors and a great deal is known about its structure and
properties. Early work on the value function
includes~\cite{blair1977value,blair1979value,blair1984constructive,blair1995closed}, while
the material here is based on the work
in~\cite{GuzRal07,Guzelsoy2009,HasRal14-1,HasRal14}. 

As a starting point,
consider an instance of~\eqref{eqn:second-stage-problem} with fixed
right-hand side $b$ and let us interpret the values of the variables as specifying
a numerical ``level of engagement'' in certain activities in an economic
market.
Further, let us interpret the constraints as corresponding to limitations
imposed on these activities due to available levels of certain scarce resources
(it is most natural to think of ``$\leq$'' constraints in this
interpretation).
In each row $j$ of the constraint matrix, the coefficient $G^2_{ij}$
associated with activity (variable) $i$ can then be thought of as representing
the rate at which resource $j$ is consumed by engagement in activity $i$.
In this interpretation, the optimal primal solution then specifies the level
of each activity in which one should engage in order to maximize profits (it
is most natural here to think in terms of maximization), given the fixed level
of resources $b$.

Assuming that additional resources were available, how much should one be
willing to pay? The intuitive answer is that one should be willing to pay at
most the marginal amount by which profits would increase if more of a
particular resource were made available. Mathematically, this information can
be extracted from the \emph{value function} $\phi: \Re^{m_2} \rightarrow \Re
\cup \{\pm \infty\}$ associated with~\eqref{eqn:second-stage-problem}, defined
by
\begin{equation} \label{eqn:value-function} \tag{2SVF}
\phi(\beta) =  \inf_{y \in \P_2(\beta) \cap Y} d^2 y,
\end{equation}
for $\beta \in \Re^{m_2}$.
Since this function returns the optimal profit for any given basket of
resources, its gradient at $b$ (assuming $\phi$ is differentiable at $b$) tells
us what the marginal change in profit would be if the level of resources
available changed in some particular direction. Thus, the gradient specifies a
``price'' on that basket of additional resources.

The reader familiar with the theory of duality for linear optimization
problems should recognize that the solution to the usual dual problem
associated with an LP of the form~\eqref{eqn:second-stage-problem} (i.e., assuming $r_2
= 0$) provides exactly this same information. In fact, we describe below that
the set of optimal solutions to the LP dual are precisely the subgradients of
its associated value function. This dual solution can hence be interpreted as
a linear price on the resources and is sometimes referred to as a vector of
``dual prices.'' The optimal dual prices allow us to easily determine whether
it will be profitable to enter into a particular activity $i$ by comparing the
profit $d^2_i$ obtained by entering into that activity to the cost $uG^2_i$ of
the required resources, where $u$ is a given vector of dual prices and $G_i^2$
is the $i$-th column of $G^2$. The difference $d^2_i - uG^2_i$ between the
profit and the cost is the \emph{reduced profit/cost} in linear optimization.
It is easily proven that the reduced profit of each activity entered into at a
non-zero level (i.e., reduced profits of the variables with non-zero value)
must be non-negative (again, in the case of maximization) and duality provides
an intuitive economic interpretation of this result.

Although the construction of the full value function is challenging even in the
simplest case of a linear optimization problem, approximations to the value
function in the local area around $b$ can still be used for sensitivity
analysis and in optimality conditions. The general dual problem we describe
next formalizes this idea by formulating the problem of constructing a
function that bounds the value function from below everywhere but yields a
strong approximation near a fixed right-hand side~$b$. Such a so-called ``dual
function'' can yield approximate ``prices'' and its iterative construction
can also be used in a more technical way to guide the evolution of an
algorithm by providing gradient information helpful in finding the optimal
solution, as well as providing a proof of its optimality.

\subsection{Dual Functions}

The above discussion leads to the natural concept of a \emph{dual (price)
function} from which we can derive a general notion of a \emph{dual
  problem}.
\begin{definition} \label{def:dual-function}
  A \emph{dual function} $F: \Re^{m_2} \rightarrow \Re$ is a function that
  satisfies $F(\beta) \leq \phi(\beta)$ for all $\beta \in \Re^{m_2}$. 
  We call such a dual function \emph{strong} at $b \in \Re^{m_2}$ if
  $F(b) = \phi(b)$.
\end{definition}
\noindent Dual functions are naturally associated with \emph{relaxations} of
the original problem, as the value function of any relaxation yields a
feasible dual function. In particular, the value function of the well-known LP
relaxation is the best convex under-estimator of the value function.

Also of interest are functions that bound the value function from above, which
we refer to as \emph{primal functions}. 
\begin{definition} \label{def:primal-function}
  A \emph{primal function} $H: \Re^{m_2} \rightarrow \Re$ is a function that
  satisfies $H(\beta) \geq \phi(\beta)$ for all $\beta \in \Re^{m_2}$. 
  We call such a primal function \emph{strong} at $b$ if
  $H(b) = \phi(b)$.
\end{definition}
\noindent In contrast to dual functions, primal functions are naturally
associated with \emph{restrictions} of the original problem and the value
function of any such restriction yields a valid primal function. 

It is immediately evident that a pair of primal and dual functions yields
optimality conditions. If we have a primal function $H^*$ and a dual function
$F^*$ such that $F^*(b) = \gamma = H^*(b)$ for some $b \in \Re^{m_2}$, then
we must also have $\phi(b) = \gamma$. Proofs of optimality of this nature are
produced by many optimization algorithms.

\subsection{Dual Problems}

The concepts we have discussed so far further lead us to the definition of a
generalized dual problem, originally introduced in~\cite{GuzRal07}, for an
instance of~\eqref{eqn:second-stage-problem} with right-hand side $b \in
\Re^{m_2}$. This problem simply calls for the construction of a dual function
that is strong for a particular fixed right-hand side $b \in \Re^{m_2}$ by
determining 
\begin{equation} \label{eqn:MILPD} \tag{MILPD} 
  \max_{F \in \Upsilon^{m_2}} \{F(b): F(\beta) \leq \phi(\beta),\ \forall
  \beta \in \Re^{m_2}\}, 
\end{equation}
where $\Upsilon^{m_2} \subseteq \{f \;|\; f : \Re^{m_2} \ra \Re\}$. Here,
$\Upsilon^{m_2}$ can be taken to be a specific class of functions, such as
linear or subadditive, to obtain specialized dual problems for particular
classes of optimization problems. It is clear that \eqref{eqn:MILPD} always
has a solution $F^*$ that is strong, provided that the value function is
real-valued everywhere (and hence belongs to $\Upsilon^{m_2}$, however it is
defined), since $\phi$ itself is a solution whenever it is finite
everywhere.\footnote{When the value function is not real-valued everywhere, we
  have to show that there is a real-valued function that coincides with the
  value function when it is real-valued and is itself real-valued everywhere
  else, but is still a feasible dual function (see~\cite{wolsey1981integer}).}

Although it may not be obvious, this notion of a dual problem naturally
generalizes existing notions for particular problem classes. For example,
consider again a parametric family of LPs defined as
in~\eqref{eqn:second-stage-problem} (i.e., assuming $r_2 = 0$). We show
informally that the usual LP dual problem with respect to a fixed instance
with right-hand side $b$ can be derived by taking $\Upsilon^{m_2}$ to be the
set of all non-decreasing linear functions in~\eqref{eqn:MILPD} and
simplifying the resulting formulation.
First, let a non-decreasing linear function $F: \Re^{m_2} \ra \Re$ be given.
Then, $\exists u \in \Re^{m_2}_+$ such that $F(\beta) = u\beta$ for
all $\beta \in \Re^{m_2}$. It follows that
\begin{align*}
  F(\beta) = u \beta \leq uG^2y = \sum_{j = 1}^{n_2} uG^2_j y_j \; \forall \beta
  \in \Re^{m_2}. 
\end{align*} 
From the above, it then follows that, for any $\beta \in \Re^{n_2}$, we have
\begin{align*}
  uG^2_j \leq d^2_j \; \forall j \in \{1, \dots, m_2\}
  & \Rightarrow u \beta \leq uG^2y \leq d^2y \; \forall y \in \P_2(\beta) \cap Y \\
  & \Rightarrow u \beta \leq \min_{y \in \P_2(\beta) \cap Y} d^2 y \\
 & \Rightarrow u \beta \leq \phi(\beta) \\
 & \Rightarrow F(\beta) \leq \phi(\beta).
\end{align*}
The conditions on the left-hand side above are exactly the feasibility
conditions for the usual LP dual and the final condition on the right is the
feasibility condition for~\eqref{eqn:MILPD}. Hence, the usual dual feasibility
conditions ensure that $u$ defines a linear function that bounds the value function
from below and is a dual function in the sense we have defined. The fact that
the epigraph of $\phi$ is a convex polyhedral cone in this case (it is the max
of linear functions associated with extreme points of the feasible region of
the dual problem) is enough to show that the dual~\eqref{eqn:MILPD} is strong
in the LP case, even when we take $\Upsilon^{m_2}$ to be the set of
(non-decreasing) linear functions. Furthermore, it is easy to show that any
subgradient of $\phi$ at $b$ is an optimal solution (and in fact, the set of
all dual feasible solutions is precisely the subdifferential of the value
function at the origin).

The concepts just discussed can be easily seen in the following small
example (note that this example is equality-constrained, in which case most of
the above derivation carries through unchanged, but the dual function no
longer needs to be non-decreasing).

\begin{example} \label{ex:lp-value-function}
\begin{align*}
	\min &\ 6y_1 + 7 y_2 + 5y_3 \\ 
	\st &\  2y_1 -7 y_2 + y_3 = b \\
	&\ y_1,y_2,y_3 \in \Re_+.
\end{align*}
The solution to the dual of this LP is unique whenever $b$ is non-zero and can
be easily obtained by considering the ratios $c_j/a_j$ of objective
coefficient to constraint coefficient for $j = 1,2, 3$, which determine which
single primal variable will take a non-zero value in the optimal basic
feasible solution. Depending on the sign of $b$, we obtain one of two possible
dual solutions:
\begin{equation*}
  u^* =  \begin{cases}
    \displaystyle 6/2 = 3 & \textrm{if } b > 0 \\
    \displaystyle 7/(-7) = -1 & \textrm{if } b < 0.
    \end{cases}
  \end{equation*}
Thus, the value function associated with this linear optimization problem is as shown in
Figure~\ref{pic:lp-value-function}. Note that, when $b = 0$, the dual solution
is not unique and can take any value between $-1$ and $3$. This set of
solutions corresponds to the set of subgradients at the single point of
non-differentiability of the value function. 
  \begin{figure}[t]
\begin{center}
  \includegraphics[scale=0.8]{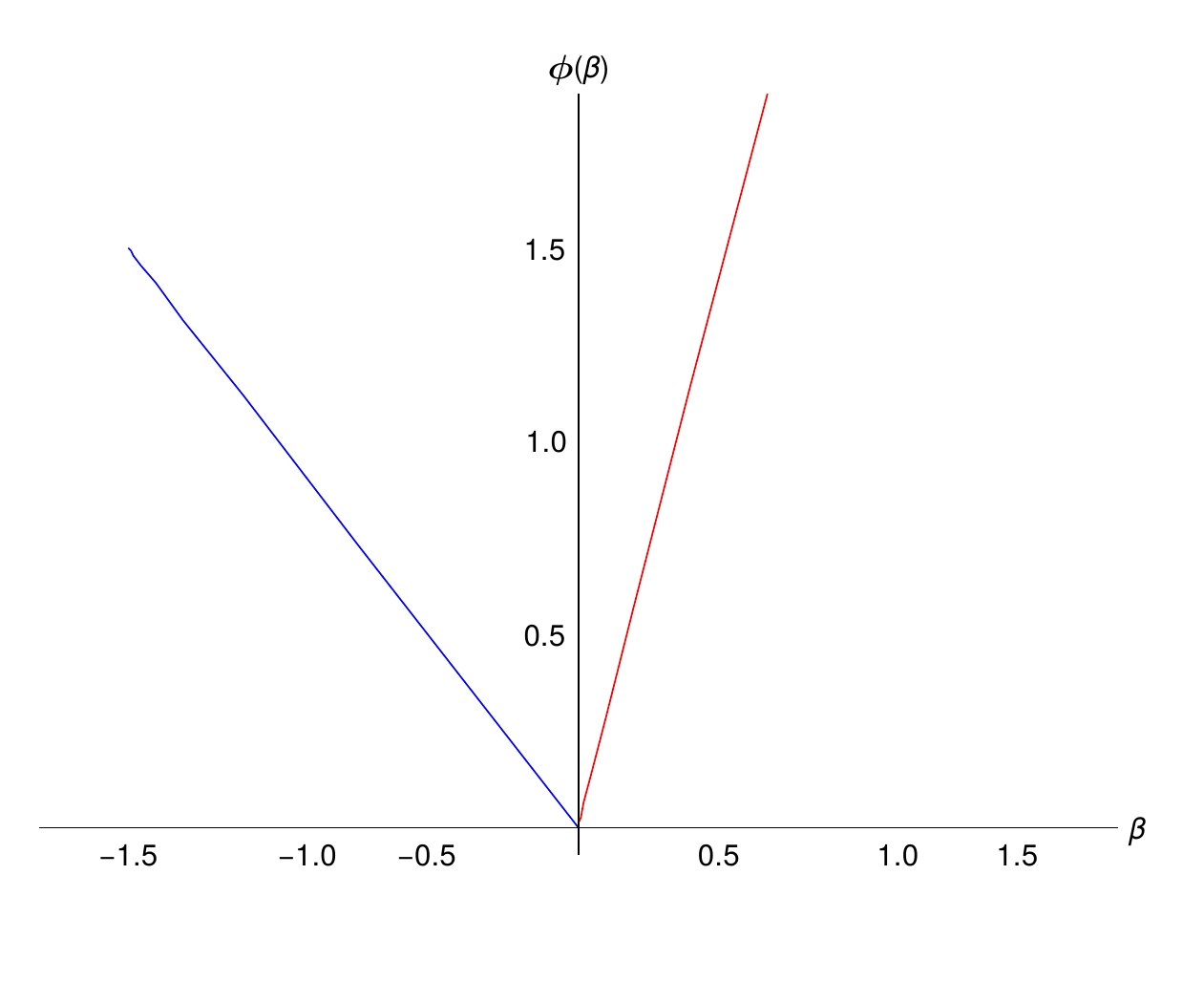}
\caption{Value Function of the LP in
  Example~\ref{ex:lp-value-function}.~\label{pic:lp-value-function}} 
\end{center}
\end{figure}
This function has one of two gradients for all points that
  are differentiable, and these gradients are equal to one of the two
  dual solutions derived above.   
\end{example}

\subsection{Structure of the Value Function \label{sec:value-function}}

As we mentioned previously, solution methods for~\eqref{eqn:2SMILP} inherently
involve the explicit or implicit approximation of several functions, including
the value function~$\phi$ in~\eqref{eqn:value-function} and the risk function
$\Xi$ in~\eqref{eqn:risk-function}, which ultimately derives its structure
from $\phi$.
Here, we summarize the results described in the series of
papers~\cite{GuzRal07,Guzelsoy2009,HasRal14-1}. Most importantly, the function
is piecewise polyhedral, lower semi-continuous, subadditive, and has a
discrete structure that is derivative of the structure of two related value
functions which we now introduce.

Let $y_C$, $d^2_C$, and $G^2_C$ be the parts of each of the vectors/matrices
describing the second-stage problem~\eqref{eqn:second-stage-problem} that are
associated with the continuous variables 
and let $y_I$, $d^2_I$, and $G^2_I$ be likewise for the integer variables.
The \emph{continuous restriction} (CR) is the LP obtained by
dropping the integer variables in the second-stage problem (or equivalently,
setting them to zero).
This problem has its own value function, defined as
\begin{equation} \label{eqn:zc} \tag{CRVF}
\phi_C(\beta) = \min_{y_C \in \Re^{n_2-r_2}_+} \{d^2_C y_C \mid G^2_C y_C \geq
\beta\}. 
\end{equation}
\noindent On the other hand, if we instead
drop the continuous variables from the problem, we can then consider the
\emph{integer restriction} (IR), which has value function 
\begin{equation} \label{eqn:zi} \tag{IRVF}
\phi_I(\beta) = \min_{y_I \in \Z^{r_2}_+} \{d^2_I y_I \mid G^2_I y_I \geq
\beta\}.  
\end{equation}
To illustrate how these two functions combine to yield the structure of $\phi$
and to briefly summarize some of the important results from the study of this
function carried out in the aforementioned papers, consider the following
simple example of a two-stage stochastic mixed integer linear optimization problem with a
single constraint. 
Note that, in this example, $d^1 = d^2$ and we are again considering the
equality-constrained case in order to make the example a bit more interesting.

\begin{example} \label{ex:milp-value-function}
\begin{equation*} 
  \begin{aligned} 
    \min \quad & \Psi(x_1, x_2) = \ -3x_1 - 4 x_2 +  \E [ \phi(b^2_\omega -  2x_1 - 0.5 x_2)]\\ 
    \st \quad
    &\ x_1 \leq 5, x_2 \leq 5\\
    &\ x_1, x_2 \in \Re_+,
  \end{aligned}
\end{equation*}
\noindent where
\begin{equation*} 
\begin{aligned}
  \phi(\beta) = \min \quad &\ 6y_1 + 4 y_2 + 3 y_3 + 4  y_4 + 5 y_5 + 7 y_6\\
  \st \quad
  &\ 2y_1 + 5 y_2 - 2y_3 - 2 y_4 + 5 y_5  + 5y_6= \beta\\
  &\ y_1,y_2,y_3 \in \Z_+,  y_4, y_5,y_6 \in \Re_+,
\end{aligned}
\end{equation*}
with $\Omega = \{1, 2\}$, $b_1^2 = 6$, and $b_2^2 = 12$.
Figure~\ref{fig:milp-value-function} shows the objective
function~$\Psi$ and the second-stage value function $\phi$ for
this example.
\begin{figure} [h]
  \centering 
  \hfill
    \includegraphics[width=0.46\textwidth]{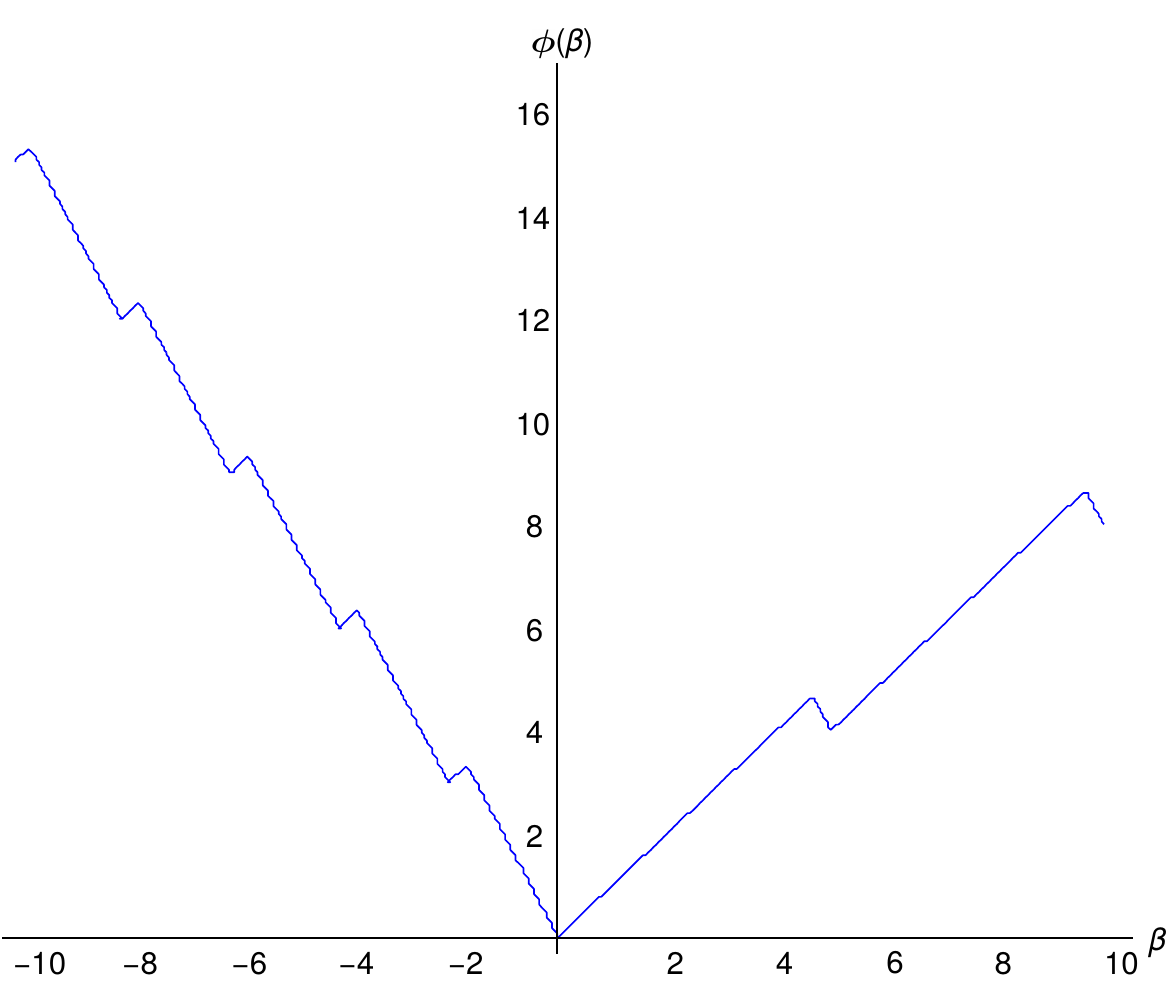}
    \includegraphics[width=0.5\textwidth]{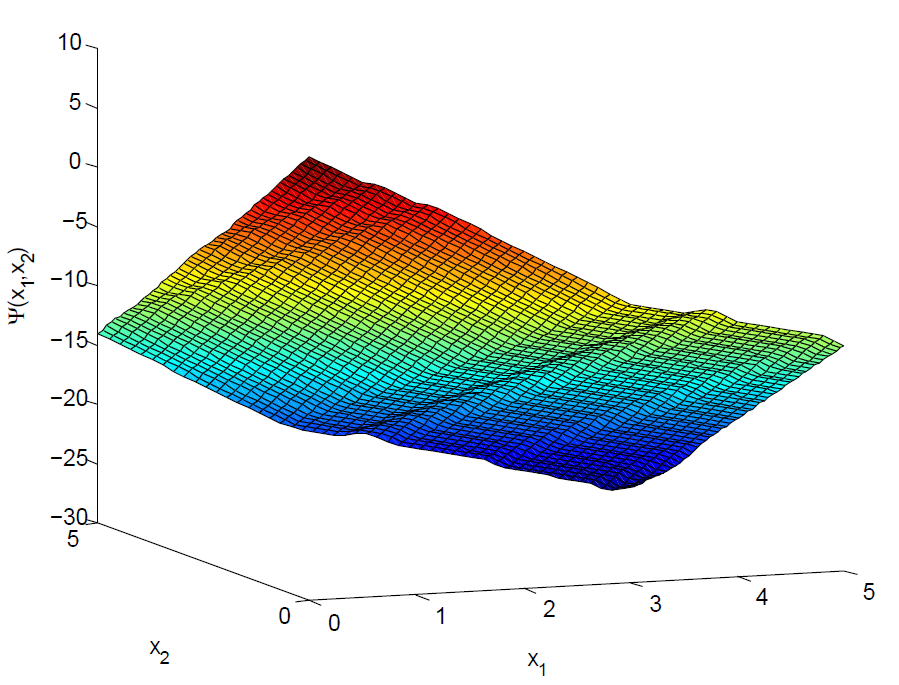}
  \hfill
  \caption{Illustration of the functions from
    Example~\ref{ex:milp-value-function}, with the second-stage value function $\phi(\beta)$ (left) and the objective function $\Psi(x_1,x_2)$ (right). \label{fig:milp-value-function}}
\end{figure}

Examining Figure~\ref{ex:milp-value-function}, it appears that $\phi$ is
comprised of a collection of translations of $\phi_C$, each of which has a structure
that is similar to the value function of the LP in
Example~\ref{ex:lp-value-function}. At points where $\phi$ is differentiable,
the gradient always corresponds to one of the two solutions to the dual of the
continuous restriction (precisely as in Example~\ref{ex:lp-value-function},
dual solutions are the ratios of objective function coefficients to constraint coefficients
for the continuous variables), which are in turn also the gradients of
$\phi_C$.
The so-called points of \emph{strict local convexity} are the points of
non-differentiability that are the extreme points of the epigraphs of the
translations of $\phi_C$ and are determined by the solutions to the integer
restriction. In particular, they correspond to points at which $\phi_I$ and
$\phi$ are coincident.
For instance, observe that, in the example, $\phi_I(5) = \phi(5) = 4$.
Finally, we can observe in Figure~\ref{ex:milp-value-function} how the structure of $\phi$
translates into that of $\Psi$.
\end{example}

Although the case illustrated in the example is of a single constraint, these
properties can be made rigorous and do generalize to higher dimension with
roughly the same intuition. The general principle is that the value function
of an MILP is the minimum of translations of the value function of the
continuous restriction $\phi_C$, where the points of translation (the
aforementioned points of \emph{strict local convexity}) are determined by the
value function of the integer restriction $\phi_I$.
\begin{theorem} \cite{HasRal14-1}
\label{the:bslc} 
Under the assumption that $\{\beta \in \Re^{m_2} \mid \phi_I(\beta) <
\infty\}$ is compact and $\epi(\phi_C)$ is pointed, there exists a finite set
$\Scal \subseteq Y$ such that 
\begin{equation*} 
  \phi(\beta) = \min_{y_I \in \Scal} \{d_I^2 y_I + \phi_C(\beta - G_I^2 y_I)\}. 
\end{equation*}
\end{theorem}
\noindent Under the assumptions of the theorem, this result provides a finite
description of $\phi$. 

\subsection{Approximating the Value Function}

Constructing functions that bound the value function of an MILP is an
important part of solution methods for both traditional single-stage
optimization problems and their multistage counterparts. Functions bounding
the value function from below are exactly the \emph{dual functions} defined
earlier and arise from \emph{relaxations} of the original problem (such as the
LP relaxation, for instance).
They are naturally obtained as by-products of common solution
algorithms, such as branch and bound, which itself works by iteratively
strengthening a given relaxation and produces a dual proof of optimality. 

Functions that bound the value function from above are the \emph{primal
  functions} defined earlier and can be obtained by considering the value
function of a \emph{restriction} of the original problem. 
While it is a little less obvious how to obtain such functions in general
(solution algorithms generally do not produce practically useful primal
functions), they can be obtained by taking the minimum over the value
functions of restrictions obtained by fixing the values of integer variables,
as we describe below.

Both primal and dual functions can be iteratively improved by
producing new such functions and combining them with existing ones by
taking the maximum over several bounding functions in the dual case or the
minimum over several bounding functions in the primal case.
When such functions are iteratively constructed to be strong at different
right-hand side values, such as when solving a sequence of instances with
different right-hand sides, such a technique can yield an approximation with
good fidelity across a larger part of the domain than any singly constructed
function could---this is, in fact, the principle implicitly behind the
algorithms we describe in Section~\ref{sec:dual-algorithms}. 

\subsubsection{Dual Functions from Branch and Bound \label{sec:bandb-dual-fcn}}

Dual functions can be obtained from most practical solution algorithms for
solving the MILP associated with~\eqref{eqn:value-function} with input $\beta
= b^2_\omega - A_\omega^2 x$, i.e., for computing $\phi(b^2_\omega -
A^2_\omega x)$. 
This is because their existence is (at least implicitly) the very source of
the proof of optimality produced by such algorithms. 
To illustrate, we show how a dual function can be obtained as the by-product
of the branch-and-bound algorithm.  

Branch and bound is an algorithm that searches the feasible region by
partitioning it and then recursively solving the resulting subproblems.
Implemented naively, this results in an inefficient complete enumeration, but
this potential inefficiency is avoided by utilizing upper and lower bounds
computed for each subproblem to intelligently ``prune'' the search. The
recursive partitioning process can be envisioned as a process of searching a
rooted tree, each of whose nodes corresponds to a subproblem. Although it is
not usually described this way, the branch-and-bound algorithm can be
interpreted as constructing a function feasible to~\eqref{eqn:MILPD}, thus
producing not only a solution but also a dual proof of its optimality.

To understand this interpretation, suppose we evaluate $\phi(b)$ for $b \in
\Re^{m_2}$ by solving the associated MILP using a standard
branch-and-bound algorithm with branching done on elementary (a.k.a. variable)
disjunctions of type $y_j \leq \pi_0 \vee y_j \geq \pi_0 + 1$ for some $j \in
\{1, \dots, r_2\}$ and $\pi_0 \in \mathbb{Z}$. Because the
individual subproblems in the branch-and-bound tree differ only by the
bounds on the variables, it will be convenient to assume that
all integer variables have finite initial lower and upper bounds denoted by
the vectors $u,l \in \Z^{r_2}$ (such bounds exist by Assumption 1, even if
they are not part of the formulation). In this case, we have that
\begin{equation*} 
  \P_2(\beta) = \{y \in \Re^{n_2} \mid G^2y \geq \beta,  y_I \geq l, -y_I \geq
  -u\}.
\end{equation*}
The solution of~\eqref{eqn:second-stage-problem} by branch and bound for
right-hand side $b$ yields a branch-and-bound tree whose leaves, contained in
the set $\T$, are each associated with the subproblem
\begin{equation*}
\min \{d^2 y \mid y \in \P_2^t(b) \cap Y\},
\end{equation*}
where
\begin{equation*}
  \P_2^t(\beta) = \{ y \in \Re^{n_2} \mid G^2y \geq \beta, y_I \geq l^t, -y_I \geq -u^t\} 
\end{equation*}
is the parametric form of the polytope containing the feasible points
associated with subproblem $t \in \T$, with $l^t, \ u^t \in \Z_+^{r_2}$ being
the modified bounds on the integer variables imposed by branching. 

Assuming that no pruning took place, the validity of the overall method comes
from the fact that valid methods of branching ensure that
\begin{equation*}
\bigcup_{t \in \T} \P_2^t(b) \cap Y = \P_2(b) \cap Y,
\end{equation*}
so that the feasible regions of the leaf nodes constitute a partition of the
original feasible region. This partition can be thought of as a single logical
disjunction that serves to strengthen the original LP relaxation. 
The proof of optimality that branch and bound produces derives from global
lower and upper bounds derived from local bounds associated with each node $t
\in \T$. We denote by $L^t$ and $U^t$, respectively, the lower and upper
bounds on the optimal solution value of the subproblem, where
\begin{align*}
  L^t & = \phiLP^t(b), \quad
  U^t = d^2 y^t,
\end{align*}
in which $\phi^t_{LP}$ is as defined in~\eqref{eqn:node-value-function} below
and $y^t \in Y$ is the best solution known for subproblem $t \in \T$ ($U^t =
\infty$ if no solution is known). Assuming~\eqref{eqn:second-stage-problem} is
solved to optimality and $y^* \in Y$ is an optimal solution, we must have 
\begin{align*}
  L := \min_{t \in \T} L^t = d^2 y^* = \min_{t \in \T} U^t =: U,
\end{align*}
where $L$ and $U$ are the global lower and upper bounds.

From the information encoded in the branch-and-bound tree, the overall dual
function can be constructed by deriving a parametric form of the lower bound,
combining dual functions for the individual subproblems in set $\T$.
For this purpose, we define the value function
\begin{equation}\label{eqn:parametric-node-dual} \tag{PNVF}
  \phiLP(\beta, \lambda, \nu) = \min \{d^2 y \mid y \in \P_2(\beta), \lambda
  \leq y_I \leq \nu, y_C \geq 0\}
\end{equation}
of a generic LP relaxation, which captures the bounds on the integer variables
as also being parametric. Using~\eqref{eqn:parametric-node-dual}, the value
function of the LP relaxation associated with a particular node $t \in \T$
(only parametric in the original right-hand side) can be obtained as
\begin{equation} \label{eqn:node-value-function} \tag{NVF} 
  \phiLP^t(\beta) = \phiLP(\beta, l^t, u^t).
\end{equation}
For all $t \in \T $ such that $\phiLP^t(b) <
\infty$, let $({v^t}, \underline{v}^t, \bar{v}^t)$ be an optimal solution to
the dual of the LP relaxation at node $t \in \T$, where $v^t$,
$\underline{v}^t$, and $\bar{v}^t$ are, respectively, the dual variables
associated with the original inequality constraints, the lower bounds on
integer variables, and the upper bounds on integer variables. Then, by LP
duality we have that
\begin{equation*} 
  \phiLP^t(b) = v^t b +  \underline{v}^t l^t - \bar{v}^t u^t.
\end{equation*}
For each node $t \in \T$ for which $\phiLP^t(b) = \infty$ (the associated
subproblem is infeasible), we instead let $({v^t}, \underline{v}^t,
\bar{v}^t)$ be a dual feasible solution that provides a finite bound exceeding
$U$ (such can be found by, e.g., adding some multiple of the dual ray that
proves infeasibility to the feasible dual solution found in the final
iteration of the simplex algorithm).

Finally, from the above we have that
\begin{equation} \tag{NDF} \label{eqn:node-dual-function}
  \uphiLP^t(\beta) = v^t \beta +  \underline{v}^t l^t - \bar{v}^t
  u^t
\end{equation}
is a dual function w.r.t. the LP relaxation at node $t$ that is strong at $b$.
Finally, we can combine these individual dual functions together to obtain 
\begin{equation}\label{eqn:bb-dual-function} \tag{BB-DF}
  \underline{\phi}^\T(\beta) = \min_{t \in \T} \uphiLP^t(\beta) =
  \min_{t \in \T} \{v^t \beta + \underline{v}^tl^t - \bar{v}^tu^t\},
\end{equation}
a dual function for the second-stage problem yielded by the tree $\T$ and
which is also \emph{strong} at the right-hand side~$b$, i.e.,
$\underline{\phi}^\T(b) = \phi(b)$.

In principle, a stronger
dual function can be obtained by replacing the single linear dual
function (which is strong at $b$ for the LP relaxation) associated with each
subproblem above by its full value function $\phiLP^t$ to obtain 
\begin{equation}\label{eqn:bb-dual-function-bis} \tag{BB-DF-BIS}
  \underline{\phi}^\T_*(\beta) = \min_{t \in \T} \phiLP^t(\beta).
\end{equation}
In practice, constructing a complete description of $\underline{\phi}^\T_*$ is
not practical (even evaluating it for a given $\beta$ requires the solution of
$|\T|$ LPs). We can instead construct a function that bounds \emph{it} from
below (and hence is also a dual function for the original problem) by
exploiting the entire collection of dual solutions arising during the solution
process. For example, let
\begin{align*}
\uphiLP(\beta, \lambda, \nu) = \max_{t \in \T} \{v^t \beta +
\underline{v}^t \lambda - \bar{v}^t \nu\}, 
\end{align*}
which consists of an approximation of the full value function $\phiLP$ using
the optimal dual solutions at each leaf node.
Replacing $\uphiLP^t(\beta)$ with $\uphiLP(\beta, l^t, u^t)$
in~\eqref{eqn:bb-dual-function} results in a potentially stronger but
still practical dual function. Of course, it is also possible to add dual
solutions found at non-leaf nodes, as well as other suboptimal dual solutions
arising during the solution process, but there is an obvious trade-off between
strength and tractability.
More details on this methodology are contained in~\cite{GuzRal07,HasRal14}.

\subsubsection{Iterative Refinement}
In iterative algorithms such as those we introduce in
Section~\ref{sec:algorithms}, the single dual
function~\eqref{eqn:bb-dual-function} we get by evaluating the value function
for one right-hand side can be iteratively augmented by taking the maximum
over a sequence of similarly derived dual functions. 
Taking this basic idea a step further, \cite{RalGuz05,RalGuz06,Guzelsoy2009}
developed methods of warm starting the solution process of an MILP. 
Such methods may serve to enhance tractability, though this is still an active area
of research.
When evaluating $\phi$ repeatedly for different values in its domain, we do
not need to solve each instance from scratch---it is possible to use the
tree resulting from a previous solve as a starting point and simply further
refine it to obtain a dual function that remains strong for the previous
right-hand side of interest and is made to be strong for a new right-hand
side. 

\cite{HasRal14} shows how to use this iterative-refinement approach
to construct a lower approximation of the value function of an MILP in the
context of a Benders-like algorithm for two-stage stochastic optimization
within a single branch-and-bound tree. 
In fact, with enough sampling this method can be used to construct a single
tree whose associated dual function is strong at every right-hand side
(provided the set of feasible right-hand sides is finite).
refinement approach in approximating the value function of
Example~\ref{ex:milp-value-function}. 

\begin{figure}[h!]
\begin{minipage}{0.55\textwidth}
\centering
\raisebox{-\height}{
\begin{tikzpicture}[scale=0.5]
\tikzset{level distance=90pt,
    sibling distance=1.5cm,
    level 1/.style={sibling distance=.8cm},
    level 2/.style={sibling distance=0pt},
    level 3/.style={sibling distance=0pt},
    execute at begin node=\strut,
    every tree node/.style={align=center},
    edge from parent/.append style={very thick,-stealth}
    }
\Tree
[.{Node 0 \\ $\uphiLP^0 = 0.8 \beta$}   
	\edge [DarkGreen] 	node[auto=right] {\color{black} $x_2 = 0$};
    [.{Node 1 \\ $\uphiLP^1 = \beta$}
    ]
  	\edge[DarkGreen]	node[auto=left] {\color{black} $x_2 \geq 1$};
    [.{Node 2 \\ $\uphiLP^2 = -1.5 \beta+11.5 $} 
    ] 
]
\end{tikzpicture}}
\raisebox{-\height}{
\begin{tikzpicture}[scale=0.5]
\tikzset{level distance=90pt,
    sibling distance=1.5cm,
    level 1/.style={sibling distance=.8cm},
    level 2/.style={sibling distance=1 cm},
    level 3/.style={sibling distance=0pt},
    execute at begin node=\strut,
    every tree node/.style={align=center},
    edge from parent/.append style={very thick,-stealth}
    }
\Tree
[.{Node 0 \\ $\uphiLP^0 = 0.8 \beta$}   
	\edge [DarkGreen] 	node[auto=right] {\color{black} $x_2  = 0$};
    [.{Node 1 \\ $\uphiLP^1 = \beta $}
    ]
  	\edge[DarkGreen]	node[auto=left] {\color{black} $x_2 \geq 1$};
    [.{Node 2 \\ $\uphiLP^2 = -1.5 \beta + 11.5$} 
	\edge [red] 	node[auto=right] {\color{black} $x_2  = 1$};
 	[.{Node 3 \\ $\uphiLP^3 = \beta -1 $}
         ]
	\edge [red] 	node[auto=left] {\color{black} $x_2  \geq 2$};
        [.{Node 4 \\ $\uphiLP^4 = -1.5 \beta + 23$} ]
    ] 
]
\end{tikzpicture}
}
\end{minipage}
\begin{minipage}{0.45\textwidth}
  \centering 
  \includegraphics[width =0.9 \textwidth]{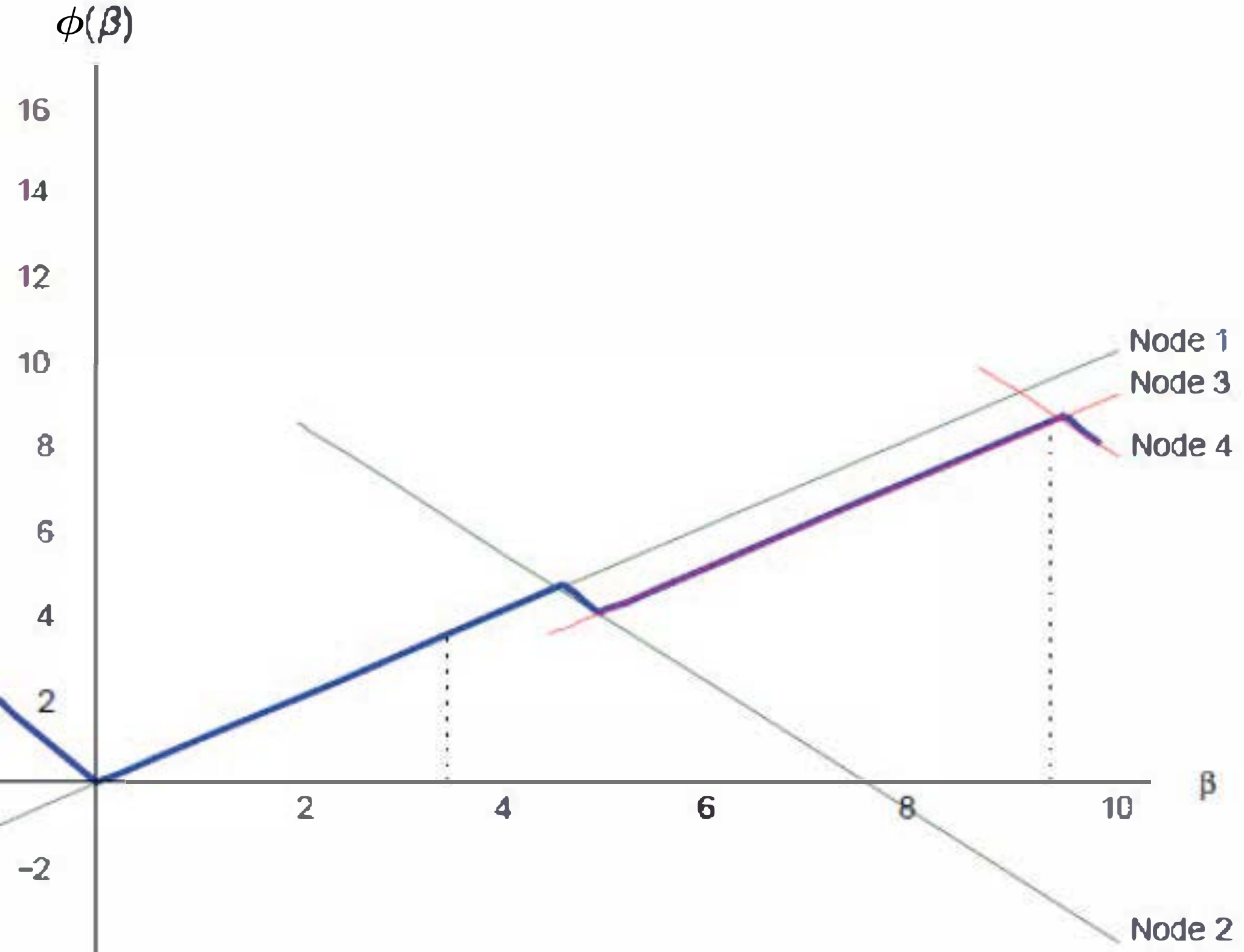}
\end{minipage}
\caption{Approximation of the value function from
  Example~\ref{ex:milp-value-function} as described in Example~\ref{Ex.BBtree}.} \label{smallBBMILP}
\end{figure}
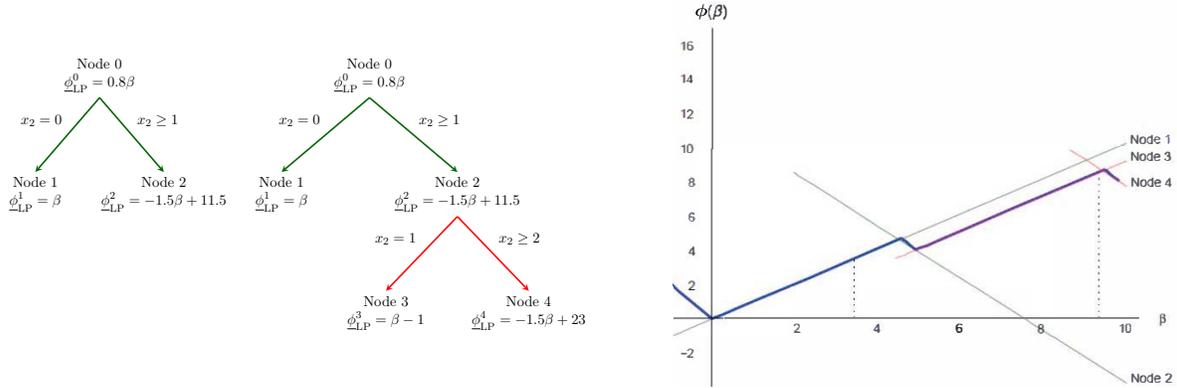

\begin{example} \label{Ex.BBtree}
  Consider the value function of Example~\ref{ex:milp-value-function}, reported in Figure~\ref{fig:milp-value-function}.
  The sequence of evaluations of the value function
  in this example are the ones arising from first-stage solutions generated by
  solving the master problem in a generalized Benders algorithm, such as the
  one described in Section~\ref{sec:dual-algorithms}.  Here, we only illustrate the way in which the dual function is iteratively refined in each step.
  
  We first evaluate $\phi(3.5)$ by branch and bound. Figure~\ref{smallBBMILP}
  shows both the tree obtained (far left) and the value function itself (in blue).
  The dual function arises from the solution to the dual of the LP
  relaxation in each of the nodes in the branch-and-bound tree.
  We exhibit the values of the dual solution for
  each node in the tree in Table~\ref{tab:nodal-duals}. Explicit upper and lower bounds were added with upper
  bounds initially taking on a large value $M$, representing infinity. Note that the dual values associated with the bound constraints are actually
  nothing more than the reduced costs associated with the variables.
  \begin{table}[h!]
    \caption{Dual solutions for each node in the branch-and-bound tree}
    \label{tab:nodal-duals}
    \setlength{\tabcolsep}{5pt}
    \begin{center}
      \begin{tabular}{r|r|rrrrrr|rrrrrr}
        $t$ & $v^t$ & \multicolumn{6}{c|}{$\underline{v}^t$} & \multicolumn{6}{c}{$\bar{v}^t$} \\
        \hline
  $0$ & $ 0.8$ & $4.4$ & $0.0$ & $4.6$ & $5.6$ & $1.0$ & $3.0$
               & $0.0$ & $0.0$ & $0.0$ & $0.0$ & $0.0$ & $0.0$ \\
  $1$ & $ 1.0$ & $4.0$ & $0.0$ & $5.0$ & $6.0$ & $0.0$ & $2.0$
               & $0.0$ & $-1.0$ & $0.0$ & $0.0$ & $0.0$ & $0.0$ \\
  $2$ & $-1.5$ & $9.0$ & $11.5$ & $0.0$ & $1.0$ & $12.5$ & $14.5$
               & $0.0$ & $0.0$ & $0.0$ & $0.0$ & $0.0$ & $0.0$ \\
      \end{tabular}
    \end{center}
  \end{table}
  
  The dual function associated with this first branch-and-bound tree is the
  minimum of the two linear functions shown in Figure~\ref{smallBBMILP} in
  green and labeled as ``Node 1'' and ``Node 2.'' Formally, this dual function
  is
\begin{equation*}
    \underline{\phi}^{\T_1} =  \min \{ \uphiLP^1,
    \uphiLP^2\},  
  \end{equation*}
  where the nodal dual functions for the three nodes are
  \begin{align*}
    \uphiLP^0(\beta) & = 0.8 \beta \\
    \uphiLP^1(\beta) &= \beta \\
    \uphiLP^2(\beta) &= -1.5\beta + 11.5.
\end{align*}
  In other words, we have $v^1_0 = 1$ (the value of the dual variable associated
  with the single equality constraint in Node 1), while $\underline{v}^1 l^1 -
  \bar{v}^1 u^1 = 0$ (this is the contribution from the dual value
  corresponding to the bound constraints, which we take to be a constant here,
  as in~\eqref{eqn:node-dual-function}). Similarly, $v^2_0 = -1.5$ and
  $\underline{v}^1 l^1 - \bar{v}^1 u^1 = 11.5$. The dual function
  $\phiLP^{\T_1}$ is strong in the interval $[0, 5]$, but yields a weaker
  lower bound outside this interval.
  If we subsequently evaluate the right-hand side $9.5$, we see that
  \begin{equation*}
    \uphiLP^{\T_1}(9.5) = \min \{9.5, -2.75\} = -2.75 \not=
    \phi(9.5) = 8.5. 
    \end{equation*}
  To obtain a strong dual function for the new right-hand side, we identify
  that node~2 is the node whose bound needs to be improved by further
  refining the tree by branching (this is the linear function yielding the
  bound in this part of the domain). By further partitioning the subproblem
  associated with node 2, we obtain the tree pictured to the right of the
  first tree in Figure~\ref{smallBBMILP}. 
  We obtain the dual function
  \begin{equation*}
    \uphiLP^{\T_2} = \min\{ \uphiLP^1,
    \uphiLP^3,  \uphiLP^4\},
  \end{equation*}
  which is strong at the right-hand side $9.5$.

  Note that this new function is no longer strong at the initial right-hand
  side of 3.5. To ensure that this single dual function remains
  strong for all previously evaluated right-hand sides, we must take the max
  over the collection of dual functions found at each iteration. This function
  is still obtained from the single tree, but we are effectively strengthening
  the leaf node dual functions by taking the max over all dual solutions
  arising on the path from the root subproblem (this is still a bound on the
  optimal solution value to the LP relaxation). In this case, we get the
  strengthened function
  \begin{equation*}
    \min\{\max\{ \uphiLP^1, \uphiLP^0\}, \max\{
    \uphiLP^3, \uphiLP^2, \uphiLP^0\}, \max\{
    \uphiLP^4, \uphiLP^2, \uphiLP^0\}\}.
  \end{equation*}
  This can be seen as an approximation of $\phi^\T_*$ by replacing the full
  value function at each node with an approximation made of just the dual
  solutions arising on the path to the root node.
\end{example}

\subsubsection{Primal Functions} \label{subusec:upper-approx}

Upper approximations of $\phi$ can be obtained by considering the value
functions of the second-stage problem~\eqref{eqn:second-stage-problem} obtained by fixing variables. 
For example, consider the {\em integer-fixing value function}
\begin{equation} \label{eqn:boundingfcn} \tag{IFVF}
\bar{\phi}_{\hat{y}}(\beta) =  d_I^2 \hat{y}_I + \phi_C(\beta - A_{\omega I}^2 \hat{y}_I)
\end{equation}
obtained by fixing the integer part $\hat{y}_I \in \Z^{r_2}$ to be equal to
that from some previously found second-stage solution $\hat{y} \in Y$, where
$\phi_C$ is as defined in~\eqref{eqn:zc}. 
Then, we have $\bar{\phi}_{\hat{y}_I}(\beta) \geq \phi(\beta)$ for all $\beta \in \Re^{m_2}$.
If $\hat{y}$ is an optimal solution to the second-stage problem with respect
to a given first-stage solution $x \in \P_1 \cap X$ and a given realized value
of $\omega$, then we have $\hat{y} \in \argmin_{y \in \P_2(b^2_\omega -
A^2_\omega x) \cap Y} d^2 y$ and $\bar{\phi}_{\hat{y}}$ is strong at $\beta =
b^2_\omega - A^2_\omega x$.

In a fashion similar to a cutting plane method, we can iteratively improve the
global upper bounding function by taking the minimum of all bounding functions
of the form \eqref{eqn:boundingfcn} found so far, i.e., 
\begin{equation*}
\bar{\phi}(\beta) = \min_{y \in \mathcal{R}} \bar{\phi}_y(\beta),
\end{equation*}
where $\mathcal{R}$ is the set of all second-stage solutions that have been found when
evaluating $\phi(\beta)$ for different values of $\beta$.  

\begin{figure}[h!]
\begin{center}
  \scalebox{0.85}{
    \input{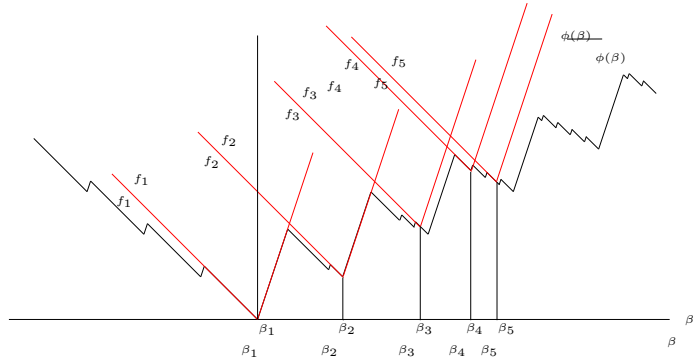}
  }
  \caption{Upper bounding functions obtained 
at right-hand sides \m{\beta_i, i=1, \dots, 5}. 
\label{fig:value_upper_apprx}}
\end{center}
\end{figure}

A pictorial example of this type of upper bounding function is shown in
Figure~\ref{fig:value_upper_apprx}, where each of the labeled cones shown is
the value function of a restriction of the original MILP. 
The upper bounding function is the minimum over all of these cones. 

\subsection{Reaction and Risk Functions}

Because it is particularly relevant in the present context, we also now
introduce a function known as the \emph{second-stage (optimistic) reaction
  function}. This function is closely related to the risk function but its
input is a second stage right-hand side, rather than a first-stage solution.
Although this function, like the second-stage value function $\phi$, takes a
right-hand side $\beta$ as its input, it can nevertheless be used to evaluate
a first-stage solution $x \in X$ in scenario $\omega \in \Omega$
by evaluating it with respect to $\beta_\omega(x)$. The function
is defined as
\begin{equation} \label{eqn:reaction-function} \tag{ReF}
\rho(\beta) = \inf \big\{ d^1 y \mid y \in\argmin \{d^2y \mid y \in
\P_2(\beta) \cap Y\} \big\},
\end{equation}
for $\beta \in \Re^{m_2}$. Note that, although the evaluation of this function
appears to require solving a bilevel optimization problem, its evaluation is
actually equivalent to solving a lexicographic optimization problem,
a somewhat easier computational task.

The importance of the function $\rho$ is to enable us to see that,
for~\eqref{eqn:2SMILP}, the scenario risk functions $\Xi_\omega$ defined
in~\eqref{eqn:second-stage-function} are not in fact completely independent
functions but, rather, are connected, since
\begin{equation*}
  \Xi_\omega(x) = \rho(\beta_\omega(x)).
\end{equation*}
Thus, these functions only differ from each other in the affine map
$\beta_\omega(x) = b^2_\omega - A_\omega^2 x$ that is applied to $x$.

The structure of the functions $\rho$ and $\Xi$ derives from that of $\phi$
and can be understood through a somewhat more involved application of the same
principles used to derive the function $\phi$.
Their structure, though
combinatorially much more complex, is nevertheless also piecewise polyhedral.
Approximations of $\Xi$ can be derived easily from approximation of $\rho$ in
a straightforward way,
since $\Xi = \E \left[\Xi_\omega\right]$ and the scenario risk functions are
themselves defined in terms of $\rho$, as discussed earlier.
The approximation of $\rho$ is quite involved, but it can be obtained
by methods that are natural generalizations of those used to approximate
$\phi$. The main challenge is that the evaluation of $\rho(\beta)$ for a
particular value of $\beta$ itself reduces to the solution of a lexicographic
optimization problem, which in turn requires knowledge of $\phi$. We may
approximate $\rho$ by repeatedly evaluating it, extracting primal and dual
information from the solution algorithm as we do with $\phi$, but this
requires repeatedly evaluating~$\phi$, which is itself expensive.

In the algorithms we discuss in Section~\ref{sec:algorithms}, we approach this
difficulty by constructing a single approximation of $\phi$ in tandem with the
approximation of $\rho$. We need only ensure that the approximation of $\phi$
is guaranteed to be strong (i.e., equal to the value function) exactly in the
region needed to properly evaluate $\rho$. The result is an approximation of
$\rho$ that is strong in the region of a given right-hand side and this is
exactly what is needed for a Benders-type algorithm to converge.
More details regarding the Benders-type algorithm are contained next in
Section~\ref{sec:algorithms}.
Further details on the structure of and methods for approximating $\rho$ and
$\Xi$ can be found in~\cite{BolRal20},
which describes a Benders-type algorithm for solving
\eqref{eqn:2SMILP}.

\subsection{Optimality Conditions}

To solidify the connection between the notion of duality described in this
section and captured in the dual problem~\eqref{eqn:MILPD}, we end this
section by formally stating both the weak and strong duality properties
arising from this theory. These properties are a proper generalization of the
well-known ones from linear optimization and can be used to derive optimality
conditions that generalize those from LP duality. These are, in turn, the
conditions that can be used to derive the reformulations presented in
Section~\ref{sec:reformulations}.
\begin{theorem}[\emph{Weak Duality}]
  If $F \in \Upsilon^{m_2}$ is feasible for~\eqref{eqn:MILPD} and $y \in
  \P_2(b) \cap Y$, then $F(b) \leq d^2 y$.
\end{theorem}
\begin{theorem}[\emph{Strong Duality}] \label{thm:strong-duality}
Let $b \in \Re^{m_2}$ be such that $\phi(b) < \infty$. Then, there exists both
$F \in \Upsilon^{m_2}$ that is feasible for~\eqref{eqn:MILPD} and $y^* \in
\P_2(b) \cap Y$ such that $F(b) = \phi(b) = d^2 y^*$.
\end{theorem}
The form of the dual~\eqref{eqn:MILPD} makes these properties rather
self-evident, but Theorem~\ref{thm:strong-duality} nevertheless yields
optimality conditions that are useful in practice. In particular, the dual
functions arising from branch-and-bound algorithms that were described earlier
in Section~\ref{sec:bandb-dual-fcn} are the strong dual functions that provide
the optimality certificates for the solutions produced by such algorithms and
are the basis on which the algorithms described in Section~\ref{sec:dual-algorithms}
are developed.

\section{Reformulations \label{sec:reformulations}}

A crucial step in deriving solution algorithms is to consider some conceptual
reformulations that underlie the problems under study, 
each of which suggests a particular algorithmic strategy. These formulations
are heavily based on the duality theory and methodology in the previous
section, as
we better clarify in the following.
In all cases, the goal is to derive, from some initial problem description, a
formulation that is, at least in principle, a single-stage mathematical
optimization problem that can be tackled with (possibly generalized versions
of) the standard algorithmic approaches used for solving mathematical
optimization problems. In this case, as in the solution of single-stage MILPs,
the main tools are cutting plane methods (branch and cut) and decomposition
methods (Benders' algorithm and the exploitation of the block structure using
a Dantzig-Wolfe decomposition).

\subsection{Value-Function (Optimality-based) Reformulation}

The first reformulation we describe is a variation
on~\eqref{eqn:2SMILP-standard}, the standard formulation used in most of the
bilevel optimization literature. This formulation
introduces the second-stage variables explicitly and formulates the problem in
the form of a classical mathematical optimization problem using a technique
that is standard---replacing the requirement that the solution to the
second-stage problem be optimal with explicit optimality conditions. To
achieve this, we introduce a constraint involving the value function $\phi$,
as well as the second-stage feasibility conditions. This is roughly equivalent
to imposing primal and dual feasibility along with equality of the primal and
dual objectives in the linear optimization case (the constraint involving the
value function must be satisfied at equality, though it is stated as an
inequality). The formulation is as follows.
\begin{mysubequations}{VFR} \label{eqn:value-function-reformulation}
\begin{align} 
\min \quad & c x + \sum_{\omega \in \Omega} p_\omega d^1 \yo \nonumber \\
\st \quad & A^1 x \geq b^1 \\
& G^2 \yo \geq b^2_\omega - A^2_\omega x \quad & &\forall \omega \in \Omega \\
& d^2 \yo \leq \phi(b^2_\omega- A^2_\omega x) \quad & & \forall \omega \in
\Omega \label{eqn:vfr-vf} \\
& x \in X \label{eqn:vfr-int-x} \\
& \yo \in Y \quad & & \forall \omega \in \Omega. \label{eqn:vfr-int-y}
\end{align}
\end{mysubequations}
It is clear that this formulation cannot be constructed explicitly, but
rather, must be solved by the iterative approximation of the constraints
involving the value function (which we refer to as the second-stage optimality
constraints). This reformulation suggests a family of methods described in
Section~\ref{sec:algorithms} in which we replace $\phi$ with a
primal function $\bar{\phi}$, as defined in
Definition~\ref{def:primal-function}).

Notice that, when $r_2 = 0$, so that the second-stage
problem~\eqref{eqn:second-stage-problem} is a linear optimization problem, we can exploit
the fact that the optimality conditions for this problem involve linear
functions. This allows for, in essence, substituting for $\phi$ the objective
function of the classical LP dual of~\eqref{eqn:second-stage-problem}, after
introducing the corresponding variables and constraints. This, overall, leads
to a tractable \emph{primal-dual} reformulation---the technique is applied,
for instance, in~\cite{coniglio2015generation}. The alternative idea of,
rather than the dual of~\eqref{eqn:second-stage-problem}, introducing its KKT
conditions, is arguably more popular and has been often exploited in a number
of ``classical'' works on mixed integer bilevel optimization problems,
including, among others,~\cite{labbe98}. Note, however, that
while there is an analog of this reformulation that applies in
the setting of~\eqref{eqn:2SMILP} (see~\cite{DeNegre2011}),
it has so far proved not to be practical and, therefore, we will not present
  any algorithms for its solution in Section~\ref{sec:algorithms}.

\subsection{Risk-Function (Projection-based) Reformulation}
  
The next reformulation we consider exploits the
finiteness of $\Omega$ and avoids introducing the second-stage variables
explicitly.
It reads as follows.
\begin{equation} \label{eqn:risk-reformulation} \tag{RFR}
  \begin{aligned}
\min \quad & c^1 x + \sum_{\omega \in \Omega}p_{\omega} z^\omega \\
\st \quad & z^\omega \geq \Xi_\omega(x) & \omega \in \Omega\\
& x\in \P_1 \cap X\\
& z^\omega \in \mathbb{R} & \omega \in \Omega.
  \end{aligned}
\end{equation}
This reformulation mirrors the original formulation implicitly adopted when we
first defined~\eqref{eqn:2SMILP}, in which the second-stage variables are not
(explicitly) present.
However, we can also interpret it as a projection onto the $X$-space of
the value-function reformulation described in the previous section.
In fact, it is not hard to see that the set $\{x \in \P_1 \mid \Xi(x) <
\infty\}$ is exactly $\mathcal{F}^1$ (the projection of the feasible region
onto the space of the first-stage variables) as defined in~\eqref{eqn:f1}.
As such, this formulation is a natural basis for a Benders-type algorithm that we
describe in Section~\ref{sec:dual-algorithms}, in which we replace $\Xi$ with
an under-estimator to obtain a \emph{master problem} which is then
iteratively improved until convergence. 

\subsection{Polyhedral (Convexification-based) Reformulation}
An apparently unrelated reformulation generalizes the notion of
\emph{convexification} used heavily in the polyhedral theory that underlies
the solution methodology of standard MILPs. 
Convexification considers the following conceptual reformulation:
\begin{equation}\label{eqn:polyhedral-formulation} \tag{POLY-R}
\begin{aligned} 
\min \quad & c x + \sum_{\omega \in \Omega} p_\omega d^1 \yo \\
\textrm{s.t.} \quad & (x, \yo) \in \conv(\F^\omega) \quad
\forall \omega \in \Omega,
\end{aligned}
\end{equation}
where $\F^\omega$ is the feasible region under scenario $\omega$, defined as
in~\eqref{eqn:F-omega}. 
Under our earlier assumptions, the convex hull of $\F^\omega$ is a polyhedron
whose extreme points are members of $\F^\omega$. Thus, due to the linearity of
the objective function, we can w.l.o.g. replace the requirement that $(x, \yo)
\in \F^\omega$ with the requirement that $(x, \yo) \in \conv(\F^\omega)$,
thereby convexifying the feasible region. 

With this reformulation, we have transformed the requirement that the
second-stage solution be optimal for the parameterized second-stage
problem~\eqref{eqn:second-stage-problem} into a requirement that the
combination of first- and second-stage solutions lie in a polyhedral feasible
region. 
This reformulation suggests a different class of algorithms based on the
dynamic generation of valid inequalities, such as those so successfully
employed in the case of MILPs. 
We describe an algorithm of such class in Section~\ref{sec:primal-algorithms}.

\subsection{Deterministic Equivalent (Decomposition-based) Reformulation}

Finally, we remark that the finiteness of $\Omega$ allows for solving the
problem via a block-angular reformulation based on the
formulation~\eqref{eqn:2SMILP-standard} presented earlier, which is in the
spirit of the so-called deterministic equivalent reformulation used in
two-stage stochastic optimization.
This renders the stochastic problem as a deterministic MIBLP, which can then
be solved via standard methods for that case with the requisite further
reformulations (of course, exploiting the block structure of the resulting
matrices). For details, see~\cite{Tahernejad2019}.

\section{Algorithmic Approaches}\label{sec:algorithms}

We now summarize
a range of methodologies that arise naturally from the reformulations of the
previous section.
Any practical method of solving~\eqref{eqn:2SMILP} must have as a fundamental step 
the evaluation of $\phi(\beta)$ for certain fixed values of $\beta \in
\mathbb{R}^{m_2}$, an operation which can be challenging in itself, since the
corresponding problem is $\NPcomplexity$-hard. From the evaluation
of $\phi$, both primal and dual information is obtained, which can be used to
build approximations. While some methods explicitly build such approximations,
other methods do it only implicitly. In all cases, information about the value
function that is built up through repeated evaluations can be exploited.

Similarly, in the dual methods that we describe below, the function $\Xi$ is
also evaluated for various values of $x \in X$ (or rather the function $\rho$)
and, similarly, approximations of this function can be built from primal and
dual information obtained during its evaluation. In order to develop
computationally tractable methods, a key aspect is to limit the number of such
function evaluations as much as possible and to exploit to the maximum extent
possible the information generated as a by-product of these single function
evaluations.

\subsection{Decomposition Methods \label{sec:dual-algorithms}}

Decomposition methods are based on generalizations of Benders' original method of
decomposing a given mathematical optimization problem by splitting the
variables into two subsets and forming a master problem by projecting out one
subset.
More concretely, we are simply solving a reformulation of
the form~\eqref{eqn:risk-reformulation}. 

\subsubsection{Continuous Problems} For illustration, we consider the simplest
case in which we have only continuous variables in both stages ($r_1 = r_2 =
0$) and $d^1 = d^2$. Since the first- and second-stage objectives are the same
in this case, the full problem is nothing more than a standard linear
optimization problem, but Benders' approach nevertheless applies when either
fixing the first stage variables results in a more tractable second-stage LP
(such as a min-cost flow problem).
In the Benders approach, we (conceptually) rewrite the LP as
\begin{equation*}
\min \left\{c x + \sum_{\omega \in \Omega} p_\omega \phi(\beta^\omega(x)) \midd
x \in \mathcal{P}_1 \right\}, 
\end{equation*}
where $\phi$ is the value function~\eqref{eqn:value-function}. Note that,
because $\Xi(x) = \sum_{\omega \in \Omega} p_\omega \phi(\beta^\omega(x))$,
this is just a simplification of the original formulation implicitly adopted
when we first defined~\eqref{eqn:2SMILP}.
As we observed earlier, the value function in the LP case is the maximum of
linear functions associated with the dual solutions. Recalling that we
can restrict the description to only the extreme points of the dual feasible
region, we can further rewrite the LP as
\begin{equation}
\min \left\{c x + \sum_{\omega \in \Omega} p_\omega z^\omega \midd \begin{array}{ll}x \in
\mathcal{P}_1\\
z^\omega \geq u\beta^\omega(x), u \in \mathcal{E}, \omega \in \Omega\\
z^\omega \in \mathbb{R}, \omega \in \Omega
\end{array}\right\},
\label{eqn:benders2-lp}\tag{LP}
\end{equation}
where $\mathcal{E}$ is the set of such extreme points of the dual of the
second-stage LP (which we assumed to be bounded and nonempty). Thus, the
linear constraints involving the variable $z^\omega$ (along with the fact that
$z^\omega$ is minimized) are precisely equivalent to requiring $z^\omega =
\phi(b^2_\omega - A^2_\omega x)$, so this reformulation is exactly the
formulation~\eqref{eqn:risk-reformulation} specialized to this case.

A straightforward solution approach is then to solve~\eqref{eqn:benders2-lp}
by a cutting plane algorithm, which results in the
classical $L$-shaped method for solving (continuous) stochastic linear
optimization problems~\cite{van1969shaped}.
From the point of view we have taken in this article, this method can be
interpreted as one that approximates the value function from below as the
maximum of the strong dual functions generated in each iteration. The strong
dual functions arise from the solutions to the dual of the second-stage
problem and yield what are typically called \emph{Benders cuts} (inequalities
of the form imposed in~\eqref{eqn:benders2-lp}). The Benders approach is then
to iteratively improve this approximation until convergence.

The case $d^1 \not= d^2$ is more complex. The epigraph of the value function
of the second-stage problem is no longer necessarily a polyhedral cone, and
the function itself is no longer necessarily convex.
Formulating the equivalent of~\eqref{eqn:benders2-lp} thus requires integer
variables. Alternative formulations using the related complementary slackness
optimality conditions are also commonly used (see~\cite{colson05}).

\subsubsection{Discrete Problems} For the case in which there
  are integer variables, the approach just described can be applied by
  simply replacing the linear strong dual functions (Benders' cuts) with
  strong under-estimators of the risk function constructed from dual functions
  arising from solution algorithms for the second-stage problem, such as those
  based on branch and bound described in
  Section~\ref{sec:bandb-dual-fcn}. In this approach, we work directly
  with the reformulation~\eqref{eqn:risk-reformulation}, employing the
  generalized Benders-type algorithm summarized in Figure~\ref{fig:gen-benders}
  and a master problem defined as follows.
  \begin{equation} \label{eqn:benders-master} \tag{MASTER}
    \begin{aligned}
      \min \quad & c^1 x + \sum_{\omega \in \Omega}p_{\omega} z^\omega \\
      \st \quad & z^\omega \geq \underline{\Xi}_\omega(x) & \omega \in \Omega\\
      & x\in \P_1\\
      & z^\omega \in \mathbb{R} & \omega \in \Omega.
    \end{aligned}
  \end{equation}

\begin{figure}[h!]
\begin{boxedminipage}{0.99\textwidth}
\textbf{Step 0. Initialize} $k \leftarrow 1$, $\underline{\Xi}^0_\omega(x)
  = -\infty$ for all $x \in \Q^{n_1}$, $\omega \in \Omega$. \\[0.1in]
\textbf{Step 1. Solve the Master Problem}
\begin{itemize}
\item[a)] Set $\displaystyle \underline{\Xi}_\omega = \max_{i = 0,\dots,k-1}
  \underline{\Xi}^i_\omega$ for $\omega \in \Omega$.
\item[b)] Solve~\eqref{eqn:benders-master} to obtain an optimal solution
$(x^{k}, \{z^{\omega}_k\}_{\{\omega \in \Omega\}})$.
\end{itemize}
\textbf{Step 2. Solve the Subproblem}
\begin{itemize}
\item[a)] Evaluate $\Xi_\omega(x^k)$ to obtain an optimal solution
  $y^{\omega,k}$ for $\omega \in \Omega$ and the strong under-estimator
  $\underline{\Xi}^k_\omega$. 
\item[b)] Termination check: Is $z^\omega_k = d^1 y^{\omega,k}$ for
  $\omega \in \Omega$? 
\begin{enumerate}
\item If yes, STOP. $x^k$ is an optimal solution
to \eqref{eqn:risk-reformulation}. 
\item If no, set $k \leftarrow k+1$ and go to Step 1.
\end{enumerate}
\end{itemize}
\end{boxedminipage}
\caption{Generalized Benders Algorithm for solving 2SMILPs.
  \label{fig:gen-benders}}
\end{figure}

When $d^1 = d^2$, the approximation of the scenario risk function and of the
risk function itself reduces to the direct approximation of the second-stage
value function,
and the algorithm can
be described rather succinctly. A basic version was originally proposed as the
integer $L$-shaped algorithm for two-stage stochastic optimization problems with
integer recourse by~\cite{laporte1993integer,caroe1998shaped}. The version
based on dual functions from branch and bound that we describe here is
described in~\cite{HasRal14}.

To briefly summarize, as in the LP case, we rewrite~\eqref{eqn:2SMILP} as
\begin{equation*}
  \min \left\{c x + \sum_{\omega \in \Omega} p_\omega z^\omega \midd
    \begin{array}{ll}
      x \in \mathcal{P}_1 \cap X\\
      z^\omega \geq \phi(\beta^\omega(x)) & \omega \in \Omega\\
      z^\omega \in \mathbb{R} & \omega \in \Omega
    \end{array}
      \right\}. 
\end{equation*}
By replacing $\phi$ with the maximum of a set $\mathcal{G}^\omega$ of dual functions
associated with scenario $\omega \in \Omega$ (alternatively,
we can employ one universal set of dual functions, as indicated
in~\eqref{eqn:benders2-lp} above), we obtain a
convergent Benders-like algorithm based on iteratively solving a master
problem of the form
\begin{equation*}
  \min \left\{c x + \sum_{\omega \in \Omega} p_\omega z^\omega \midd
    \begin{array}{ll}x \in
      \mathcal{P}_1 \cap X\\
      z^\omega \geq F(\beta^\omega(x)) & F \in \mathcal{G}^\omega, \omega \in \Omega\\
      z^\omega \in \mathbb{R} & \omega \in \Omega
    \end{array}\right\},
  \label{eqn:benders2-ip}
\end{equation*}
which generalizes~\eqref{eqn:benders2-lp}.
The key to making this approach work in practice is that
the dual functions we need be easily available as a by-product of evaluating
the second-stage value function $\phi$ for a fixed value of $\beta$. 

The most general case in which $d^1 \not= d^2$ is conceptually
no more complex than that described above, but the details
of the methodology for approximating $\Xi$ and in linearizing the master
problem are quite involved. The reader is referred to~\cite{BolRal20} for the details.

\subsection{Convexification-based Methods \label{sec:primal-algorithms}}
Primal algorithms are based on the implicit solution
of~\eqref{eqn:polyhedral-formulation} and generalize the well-known framework
of branch and cut that has been so successful in the MILP case. This class of
algorithms is based on the iterative approximation of $\conv(\F^\omega)$
beginning with the approximation $\P^\omega$, the feasible region in scenario
$\omega$ of the following relaxation obtained by dropping both the value-function
constraint~\eqref{eqn:vfr-vf} and the integrality
requirements~\eqref{eqn:vfr-int-x} and~\eqref{eqn:vfr-int-y}
from~\eqref{eqn:value-function-reformulation}. 
\begin{mysubequations}{LPR} \label{eqn:lp-relaxation} 
\begin{align} 
\min \quad & c x + \sum_{\omega \in \Omega} p_\omega d^1 \yo \nonumber \\
\st \quad & A^1 x \geq b^1 \\
& G^2 \yo \geq b^2_\omega - A^2_\omega x \quad & &\forall \omega \in \Omega \label{eqn:lpr-ss} \\
& x \in \Re^{n_1}_+ \\
& \yo \in \Re^{n_2}_+ \quad & & \forall \omega \in \Omega. 
\end{align}
\end{mysubequations}
Being an LP, this relaxation is easily solved, but it is not difficult to see,
however, that it is rather weak (see, e.g., Example~\ref{ex:valid-inequalities}).
A straightforward way to strengthen it is simply by including the integrality
constraints~\eqref{eqn:vfr-int-x} and~\eqref{eqn:vfr-int-y}
from~\eqref{eqn:value-function-reformulation}. 
This leads to an MILP relaxation with feasible set
\begin{align*}
  \mathcal{S}^\omega = \P^\omega \cap (X \times Y)
\end{align*}
in scenario $\omega \in \Omega$, which, while stronger, is clearly more
difficult to solve and also still potentially weak---whether adopting it is a
computationally good idea is a purely empirical question.
When the number of scenarios is large, dropping constraints~\eqref{eqn:lpr-ss}
from the relaxation may also be advantageous, since this may reduce the size
of the relaxation.

As in cutting plane methods for MILPs, the idea is to improve this initial
formulation with the addition of inequalities valid for
$\mathcal{S}^\omega$, $\F^\omega$, $\bigcup_{\omega \in \Omega} \F^\omega$, or
even $\F^1$.
In some cases, inequalities may first be derived with respect to
$\mathcal{S}^\omega$ or $\F^\omega$ for some particular scenario $\omega \in
\Omega$ and then lifted to become valid for a larger set.
Inequalities valid for $\S^\omega$ (which can be referred to as
\emph{feasibility cuts}) can be generated using any number of well-known
procedures associated with cutting plane algorithms for mixed integer
linear programming. 
Inequalities valid for $\F^\omega$ (which can be referred to as
\emph{optimality cuts}) are the more interesting case because they can
incorporate information about the value function in order to eliminate members
of $\P^\omega$ that are not two-stage feasible. 

In early work on these methods, the authors of~\cite{DeNRal09} developed
inequalities valid for $\F^\omega$ in the case for which $\Omega$ is a
singleton and the variables must all be integer ($r_1 = n_1$ and $r_2 = n_2$),
which illustrate the basic principles. When the input data are integer, a very
simple argument can be used to generate an inequality valid for $\F^\omega$
but violated by $(\hat{x}, \hat{y})$, an extreme point of $\P^\omega$ not in
$\F^\omega$, by taking advantage of the discrete nature of the feasible set.
Assuming the solution of the LP relaxation is an extreme point of $\P^\omega$,
there is thus a hyperplane supporting $\P^\omega$ and inducing a face of
dimension $0$. As such, there exist
$f \in \Re^{n_1}$, $g \in \Re^{n_2}$, and $\gamma \in \Re$ 
such that the hyperplane $\{(x, y) \in \Re^{n_1+n_2} \mid f x + g y = \gamma\}$ intersects 
$\P^\omega$ in the unique point $(\hat{x},\hat{y})$.
Hence, we have that $f x + g y \leq \gamma$ for all $(x, y) \in \P^\omega$.
Finally, since the face of $\P^\omega$ induced by this inequality does not contain
any other members of $\mathcal{S}^\omega$, we can ``push'' the hyperplane
until it meets the next integer point without separating any additional
members of $\F^\omega$. 
Hence,
\begin{equation*}
f x + g y \leq \gamma -1 
\end{equation*}
is valid for $\F^\omega$.
This procedure is similar in spirit to the Gomory procedure for standard MILPs.
It is used, for instance, in~\cite{Dempe2017}.
We next describe the method with a brief example.
\begin{example} \label{ex:valid-inequalities}
Consider the instance
\begin{equation*}
\max_x \min_y \left\{y\mid -x + y \leq 2, -2x -y \leq -2, 3x-y\leq 3, y\leq 3,
x,y\in \mathbb{Z}_+\right\},
\end{equation*}
with $|\Omega| = 1$, whose feasible region is the set $\F = \{(0, 2),
(1, 0), (2, 3)\}$ shown in Figure~\ref{fig:bcCut}. Solving the LP relaxation
yields the point $(1,3)$, which is not feasible. This point is eliminated by
the addition of the inequality $x - 2y \geq -4$, which is valid for the
feasible region $\F$ and is obtained as a strengthening of the
inequality $x - 2y \geq -5$, which is valid for the LP relaxation itself.
\begin{figure}[h!]
\begin{center}
\scalebox{0.5}{
\input{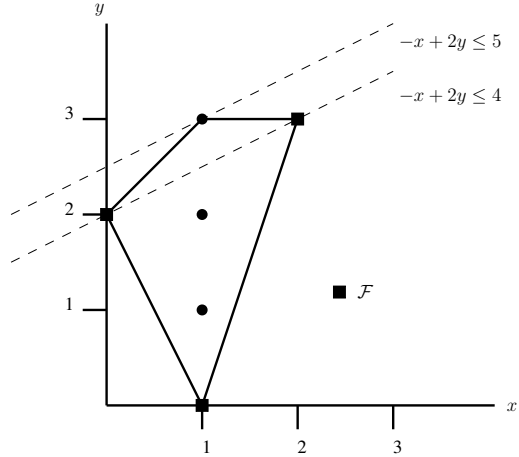}
}
\end{center}
\caption{Example of a valid inequality. \label{fig:bcCut}}
\end{figure}
\end{example}
This cut generation procedure is enough to yield a converging algorithm in this
case, but it amounts to removing infeasible points one by one and is not
scalable in general. 
An important observation is that this cut only exploits integrality of the
solutions and does not take into account any information about the second-stage
value function. 

A generalized version of this basic methodology has since been described
in~\cite{TahRalDeN16} and enhanced with additional classes of inequalities,
including those valid for the general mixed integer case described
in~\cite{fischettietal17a,fischettietal17b}.
Inequalities valid for more general discrete probability spaces are derived
in~\cite{gade2012decomposition} for the case $d^1 = d^2$. 

Stronger cuts can be obtained by using disjunctive arguments based on
knowledge of the value function.
In particular, an option is to add constraints of the form
\begin{equation*}
  d^2 y^\omega \leq \bar{\phi}(b^2_\omega - A^2_\omega x),
\end{equation*}
where $\bar{\phi}$ is a primal function, as defined in
Definition~\ref{def:primal-function}).
Such primal functions can take many forms and imposing such constraints may be
expensive. 
In general, the form of such functions will be either affine or
piecewise polyhedral (``standard'' disjunctive programming techniques can be
used to obtain a reformulation which only involves linear functions). 

\section{Conclusions}

We have introduced a unified framework for multistage mixed integer linear optimization problems which encompasses both multilevel mixed integer linear optimization problems and multistage mixed integer linear optimization problems with recourse.
Focusing on the two-stage case, we have investigated the nature of the value function of the second-stage problem and highlighted its connection to dual functions and the theory of duality for mixed integer linear optimization problems.
We have summarized different reformulations for this broad class of problems, which rely on either the risk function, the value function, or on their polyhedral nature.
We have then presented the main solution techniques for problems of this class, considering both  dual- and primal-type methods, the former based on a Benders-like decomposition to approximate either the risk function or the value function, and the latter based on a cutting plane technique that relies on the polyhedral nature of these problems.
While much work is still to be done for solving multistage mixed integer linear optimization problems with techniques that are (\emph{mutatis mutandis}, given their intrinsic harder nature) of comparable efficiency to those for solving single-level problems, we believe that the theoretical understanding of multistage mixed integer linear problems is now sufficiently mature to make this an achievable objective.

\section*{Acknowledgements}

This research was made possible with support from National Science
Foundation Grants CMMI-1435453, CMMI-0728011, and ACI-0102687, as well
as Office of Naval Research Grant N000141912330.

\bibliographystyle{plainnat}
\bibliography{MIBLP.bib}

\begin{thebibliography}{144}
\providecommand{\natexlab}[1]{#1}
\providecommand{\url}[1]{\texttt{#1}}
\expandafter\ifx\csname urlstyle\endcsname\relax
  \providecommand{\doi}[1]{doi: #1}\else
  \providecommand{\doi}{doi: \begingroup \urlstyle{rm}\Url}\fi

\bibitem[Amaldi et~al.(2013{\natexlab{a}})Amaldi, Capone, Coniglio, and
  Gianoli]{amaldi2013energy}
E.~Amaldi, A.~Capone, S.~Coniglio, and L.G. Gianoli.
\newblock Energy-aware traffic engineering with elastic demands and mmf
  bandwidth allocation.
\newblock In \emph{Proc of 18th IEEE Int. Workshop on Computer Aided Modeling
  and Design of Communication Links and Networks (CAMAD 2013)}, pages 169--174.
  IEEE, 2013{\natexlab{a}}.

\bibitem[Amaldi et~al.(2013{\natexlab{b}})Amaldi, Capone, Coniglio, and
  Gianoli]{amaldi2013network}
E.~Amaldi, A.~Capone, S.~Coniglio, and L.G. Gianoli.
\newblock Network optimization problems subject to max-min fair flow
  allocation.
\newblock \emph{IEEE Communications Letters}, 17\penalty0 (7):\penalty0
  1463--1466, 2013{\natexlab{b}}.

\bibitem[Amaldi et~al.(2013{\natexlab{c}})Amaldi, Coniglio, Gianoli, and
  Ileri]{amaldi2013single}
E.~Amaldi, S.~Coniglio, L.G. Gianoli, and C.U. Ileri.
\newblock On single-path network routing subject to max-min fair flow
  allocation.
\newblock \emph{Electronic Notes in Discrete Mathematics}, 41:\penalty0
  543--550, 2013{\natexlab{c}}.

\bibitem[Amaldi et~al.(2014)Amaldi, Coniglio, and Taccari]{amaldi2014maximum}
E.~Amaldi, S.~Coniglio, and L.~Taccari.
\newblock Maximum throughput network routing subject to fair flow allocation.
\newblock In \emph{Proc. of Int. Symp. on Combinatorial Optimization (ISCO
  2014)}, pages 1--12. Springer, 2014.

\bibitem[Amouzegar and Moshirvaziri(1999)]{amouzegar-moshirvaziri99}
M.A. Amouzegar and K.~Moshirvaziri.
\newblock Determining optimal pollution control policies: An application of
  bilevel programmingoa.
\newblock \emph{European Journal of Operational Research}, 119\penalty0
  (1):\penalty0 100--120, 1999.

\bibitem[An et~al.(2011)An, Pita, Shieh, Tambe, Kiekintveld, and
  Marecki]{an2011guards}
B.~An, J.~Pita, E.~Shieh, M.~Tambe, C.~Kiekintveld, and J.~Marecki.
\newblock Guards and {Protect}: Next generation applications of security games.
\newblock \emph{ACM SIGecom Exchanges}, 10\penalty0 (1):\penalty0 31--34, 2011.

\bibitem[Avenhaus et~al.(1991)Avenhaus, Okada, and
  Zamir]{avenhaus1991inspector}
R.~Avenhaus, A.~Okada, and S.~Zamir.
\newblock Inspector leadership with incomplete information.
\newblock In \emph{Game equilibrium models IV}, pages 319--361. Springer, 1991.

\bibitem[Bard(1983)]{bard83}
J.~Bard.
\newblock An algorithm for solving the general bilevel programming problem.
\newblock \emph{Mathematics of Operations Research}, 8\penalty0 (2):\penalty0
  260--272, 1983.

\bibitem[Bard and Moore(1992)]{bard92}
J.~Bard and J.T. Moore.
\newblock An algorithm for the discrete bilevel programming problem.
\newblock \emph{Naval Research Logistics}, 39\penalty0 (3):\penalty0 419--435,
  1992.

\bibitem[Bard et~al.(2000)Bard, Plummer, and Sourie]{bard-plummer-sourie00}
J.F. Bard, J.~Plummer, and J.C. Sourie.
\newblock A bilevel programming approach to determining tax credits for biofuel
  production.
\newblock \emph{European Journal of Operational Research}, 120:\penalty0
  30--46, 2000.

\bibitem[{Baringo} and {Conejo}(2012)]{baringo-conejo12}
L.~{Baringo} and A.J. {Conejo}.
\newblock Transmission and wind power investment.
\newblock \emph{IEEE Transactions on Power Systems}, 27\penalty0 (2):\penalty0
  885--893, May 2012.

\bibitem[Basilico et~al.(2016)Basilico, Coniglio, and
  Gatti]{basilico2016methods}
N.~Basilico, S.~Coniglio, and N.~Gatti.
\newblock Methods for finding leader-follower equilibria with multiple
  followers (extended abstract).
\newblock In \emph{Proc. of 2016 Int. Conf. on Autonomous Agents and Multiagent
  Systems (AAMAS 2016)}, pages 1363--1364, 2016.

\bibitem[Basilico et~al.(2017)Basilico, Coniglio, Gatti, and
  Marchesi]{basilico2017bilevel}
N.~Basilico, S.~Coniglio, N.~Gatti, and A.~Marchesi.
\newblock Bilevel programming approaches to the computation of optimistic and
  pessimistic single-leader-multi-follower equilibria.
\newblock In \emph{Proc. of 16th Int. Symp. on Experimental Algorithms (SEA
  2017)}. Schloss Dagstuhl-Leibniz-Zentrum fuer Informatik, 2017.

\bibitem[Basilico et~al.(2020)Basilico, Coniglio, Gatti, and
  Marchesi]{basilico2020bilevel}
N.~Basilico, S.~Coniglio, N.~Gatti, and A.~Marchesi.
\newblock Bilevel programming methods for computing
  single-leader-multi-follower equilibria in normal-form and polymatrix games.
\newblock \emph{EURO Journal on Computational Optimization}, 8:\penalty0 3--31,
  2020.

\bibitem[Ben-Ayed and Blair(1990)]{benayed90}
O.~Ben-Ayed and C.~Blair.
\newblock Computational difficulties of bilevel linear programming.
\newblock \emph{Operations Research}, 38:\penalty0 556--560, 1990.

\bibitem[Ben-Ayed et~al.(1992)Ben-Ayed, Blair, Boyce, and
  LeBlanc]{benayed-etal92}
O.~Ben-Ayed, C.~Blair, D.~Boyce, and L.~LeBlanc.
\newblock Construction of a real-world bilevel linear programming model of the
  highway network design problem.
\newblock \emph{Annals of Operations Research}, 34\penalty0 (1):\penalty0
  219--254, 1992.

\bibitem[Bertsekas(2017)]{bertsekas2017dynamic}
D.P. Bertsekas.
\newblock \emph{Dynamic Programming and Optimal Control}.
\newblock Athena scientific, 2017.

\bibitem[Birge and Louveaux(2011)]{birge2011introduction}
J.R. Birge and F.~Louveaux.
\newblock \emph{Introduction to stochastic programming}.
\newblock Springer Science \& Business Media, 2011.

\bibitem[Blair(1995)]{blair1995closed}
C.E. Blair.
\newblock A closed-form representation of mixed-integer program value
  functions.
\newblock \emph{Mathematical Programming}, 71\penalty0 (2):\penalty0 127--136,
  1995.

\bibitem[Blair and Jeroslow(1977)]{blair1977value}
C.E. Blair and R.G. Jeroslow.
\newblock The value function of a mixed integer program: I.
\newblock \emph{Discrete Mathematics}, 19\penalty0 (2):\penalty0 121--138,
  1977.

\bibitem[Blair and Jeroslow(1979)]{blair1979value}
C.E. Blair and R.G. Jeroslow.
\newblock The value function of a mixed integer program: Ii.
\newblock \emph{Discrete Mathematics}, 25\penalty0 (1):\penalty0 7--19, 1979.

\bibitem[Blair and Jeroslow(1984)]{blair1984constructive}
C.E. Blair and R.G. Jeroslow.
\newblock Constructive characterizations of the value-function of a
  mixed-integer program i.
\newblock \emph{Discrete Applied Mathematics}, 9\penalty0 (3):\penalty0
  217--233, 1984.

\bibitem[Bolusani and Ralphs(2020)]{BolRal20}
S.~Bolusani and T.K. Ralphs.
\newblock {A Framework for Generalized Benders' Decomposition and Its
  Application to Multilevel Optimization}.
\newblock Technical report, COR@L Laboratory Report 20T-004, Lehigh University,
  2020.
\newblock URL
  \url{http://coral.ie.lehigh.edu/~ted/files/papers/MultilevelBenders20.pdf}.

\bibitem[Bracken and McGill(1973)]{bracken-mcgill73}
J.~Bracken and J.T. McGill.
\newblock Mathematical programs with optimization problems in the constraints.
\newblock \emph{Operations Research}, 21\penalty0 (1):\penalty0 37--44, 1973.

\bibitem[Burgard et~al.(2003)Burgard, Pharkya, and
  Maranas]{burgard-pharkya-maranas03}
A.P. Burgard, P.~Pharkya, and C.D. Maranas.
\newblock {OptKnock: A Bilevel Programming Framework for Identifying Gene
  Knockout Strategies for Microbial Strain Optimization}.
\newblock \emph{Biotechnology and Bioengineering}, 84:\penalty0 647--657, 2003.

\bibitem[Calamai and Vicente(1994)]{calamai94}
P.~Calamai and L.~Vicente.
\newblock Generating quadratic bilevel programming problems.
\newblock \emph{ACM Transactions on Mathematical Software}, 20:\penalty0
  103--119, 1994.

\bibitem[Caprara et~al.(2016)Caprara, Carvalho, Lodi, and
  Woeginger]{capraraetal16}
A.~Caprara, M.~Carvalho, A.~Lodi, and G.~Woeginger.
\newblock Bilevel knapsack with interdiction constraints.
\newblock \emph{INFORMS Journal on Computing}, 28\penalty0 (2):\penalty0
  319--333, 2016.

\bibitem[Caramia and Mari(2015)]{caramiamari15}
M.~Caramia and R.~Mari.
\newblock Enhanced exact algorithms for discrete bilevel linear problems.
\newblock \emph{Optimization Letters}, 9\penalty0 (7):\penalty0 1447--1468,
  2015.

\bibitem[Car{\o}e and Tind(1998)]{caroe1998shaped}
C.C. Car{\o}e and J.~Tind.
\newblock {L-shaped decomposition of two-stage stochastic programs with integer
  recourse}.
\newblock \emph{Mathematical Programming}, 83\penalty0 (1):\penalty0 451--464,
  1998.

\bibitem[Castiglioni et~al.(2019{\natexlab{a}})Castiglioni, Marchesi, and
  Gatti]{castiglioni2019bealeader}
M.~Castiglioni, A.~Marchesi, and N.~Gatti.
\newblock Be a leader or become a follower: the strategy to commit to with
  multiple leaders.
\newblock In \emph{Proc. of 28th Int. Joint Conf. on Artificial Intelligence
  (IJCAI 2019)}, 2019{\natexlab{a}}.

\bibitem[Castiglioni et~al.(2019{\natexlab{b}})Castiglioni, Marchesi, Gatti,
  and Coniglio]{castiglioni2019leadership}
M.~Castiglioni, A.~Marchesi, N.~Gatti, and S.~Coniglio.
\newblock Leadership in singleton congestion games: What is hard and what is
  easy.
\newblock \emph{Artificial Intelligence}, 277:\penalty0 103--177,
  2019{\natexlab{b}}.

\bibitem[Celli et~al.(2019)Celli, Coniglio, and Gatti]{celli2019computing}
A.~Celli, S.~Coniglio, and N.~Gatti.
\newblock Computing optimal ex ante correlated equilibria in two-player
  sequential games.
\newblock In \emph{Proc. of 18th Int. Conf. on Autonomous Agents and MultiAgent
  Systems (AAMAS 2019)}, pages 909--917. International Foundation for
  Autonomous Agents and Multiagent Systems, 2019.

\bibitem[Celli et~al.(2020)Celli, Coniglio, and Gatti]{celli2020bayesian}
A.~Celli, S.~Coniglio, and N.~Gatti.
\newblock Private bayesian persuasion with sequential games.
\newblock In \emph{Proc. of 34th AAAI Conf. on Artificial Intelligence (AAAI
  2020)}, pages 1--8. AAAI Press, 2020.

\bibitem[Church and Scaparra(2006)]{ChSc06}
R.L. Church and M.P. Scaparra.
\newblock Protecting critical assets: The $r$-interdiction median problem with
  fortification.
\newblock \emph{Geographical Analysis}, 39\penalty0 (2):\penalty0 129--146,
  2006.

\bibitem[Church et~al.(2004)Church, Scaparra, and Middleton]{ChScMi04}
R.L. Church, M.P. Scaparra, and R.S. Middleton.
\newblock Identifying critical infrastructure: The median and covering facility
  interdiction problems.
\newblock \emph{Annals of the Association of American Geographers}, 94\penalty0
  (3):\penalty0 491--502, 2004.

\bibitem[Clark and Westerberg(1990)]{clark-westerberg90}
P.A. Clark and A.W. Westerberg.
\newblock Bilevel programming for steady-state chemical process design i.
  fundamentals and algorithms.
\newblock \emph{Computers \& Chemical Engineering}, 14\penalty0 (1):\penalty0
  87--97, 1990.

\bibitem[Colson et~al.(2005)Colson, Marcotte, and Savard]{colson05}
B.~Colson, P.~Marcotte, and G.~Savard.
\newblock Bilevel programming: A survey.
\newblock \emph{4OR: A Quarterly Journal of Operations Research}, 3\penalty0
  (2):\penalty0 87--107, 2005.

\bibitem[Coniglio and Gualandi(2017)]{coniglio2017separation}
S.~Coniglio and S.~Gualandi.
\newblock On the separation of topology-free rank inequalities for the max
  stable set problem.
\newblock In \emph{Proc. of 16th Int. Symp. on Experimental Algorithms (SEA
  2017)}. Schloss Dagstuhl-Leibniz-Zentrum fuer Informatik, 2017.

\bibitem[Coniglio and Tieves(2015)]{coniglio2015generation}
S.~Coniglio and M.~Tieves.
\newblock On the generation of cutting planes which maximize the bound
  improvement.
\newblock In \emph{Proc. of 14th Int. Symp. on Experimental Algorithms (SEA
  2015)}, pages 97--109. Springer, 2015.

\bibitem[Coniglio et~al.(2017)Coniglio, Gatti, and
  Marchesi]{coniglio2017pessimistic}
S.~Coniglio, N.~Gatti, and A.~Marchesi.
\newblock Pessimistic leader-follower equilibria with multiple followers.
\newblock In \emph{Proc. of 26th Int. Joint Conf. on Artificial Intelligence
  (IJCAI 2017)}, pages 171--177. AAAI Press, 2017.

\bibitem[Coniglio et~al.(2020{\natexlab{a}})Coniglio, Gatti, and
  Marchesi]{coniglio2020computing}
S.~Coniglio, N.~Gatti, and A.~Marchesi.
\newblock Computing a pessimistic stackelberg equilibrium with multiple
  followers: The mixed-pure case.
\newblock \emph{Algorithmica}, 82\penalty0 (5):\penalty0 1189--1238,
  2020{\natexlab{a}}.

\bibitem[Coniglio et~al.(2020{\natexlab{b}})Coniglio, Sirvent, and
  Weibelzahl]{coniglioAirline}
S.~Coniglio, M.~Sirvent, and M.~Weibelzahl.
\newblock Airport capacity extension, fleet investment, and optimal aircraft
  scheduling in a multi-level market model: Quantifying the costs of imperfect
  markets.
\newblock \emph{Submitted}, 2020{\natexlab{b}}.

\bibitem[Conitzer and Korzhyk(2011)]{conitzer2011commitment}
V.~Conitzer and D.~Korzhyk.
\newblock Commitment to correlated strategies.
\newblock In \emph{Proc. of 25th AAAI Conf. on Artificial Intelligence (AAAI
  2011)}, pages 632--637, 2011.

\bibitem[Conitzer and Sandholm(2006)]{conitzer2006computing}
V.~Conitzer and T.~Sandholm.
\newblock Computing the optimal strategy to commit to.
\newblock In \emph{Proc. of 7th ACM Conf. on Electronic Commerce (EC 2006)},
  pages 82--90, 2006.

\bibitem[Cook(1971)]{cook1971complexity}
S.A. Cook.
\newblock The complexity of theorem-proving procedures.
\newblock In \emph{Proceedings of the third annual ACM symposium on Theory of
  computing}, pages 151--158. ACM, 1971.

\bibitem[Cormican et~al.(1998)Cormican, Morton, and Wood]{cormican98}
K.J. Cormican, D.P. Morton, and R.K. Wood.
\newblock Stochastic network interdiction.
\newblock \emph{Operations Research}, 46\penalty0 (2):\penalty0 184--197, 1998.

\bibitem[Cote et~al.(2003)Cote, Marcotte, and Savard]{cote-marcotte-savard03}
J.-P. Cote, P~Marcotte, and G.~Savard.
\newblock A bilevel modelling approach to pricing and fare optimisation in the
  airline industry.
\newblock \emph{Journal of Revenue and Pricing Management}, 2\penalty0
  (1):\penalty0 23--36, 2003.

\bibitem[Dempe(2001)]{dempe01}
S.~Dempe.
\newblock Discrete bilevel optimization problems.
\newblock Technical Report D-04109, Institut fur Wirtschaftsinformatik,
  Universitat Leipzig, Leipzig, Germany, 2001.

\bibitem[Dempe and Kue(2017)]{Dempe2017}
S.~Dempe and F.~Mefo Kue.
\newblock Solving discrete linear bilevel optimization problems using the
  optimal value reformulation.
\newblock \emph{Journal of Global Optimization}, 68\penalty0 (2):\penalty0
  255--277, Jun 2017.
\newblock ISSN 1573-2916.
\newblock \doi{10.1007/s10898-016-0478-5}.
\newblock URL \url{https://doi.org/10.1007/s10898-016-0478-5}.

\bibitem[Dempe et~al.(2005)Dempe, Kalashnikov, and
  Rios-Mercado]{dempe-kalashnikov-mercado05}
S.~Dempe, V.~Kalashnikov, and R.~Z. Rios-Mercado.
\newblock Discrete bilevel programming: Application to a natural gas cash-out
  problem.
\newblock \emph{European Journal of Operational Research}, 166\penalty0
  (2):\penalty0 469--488, 2005.

\bibitem[DeNegre(2011)]{DeNegre2011}
S.~DeNegre.
\newblock \emph{{Interdiction and Discrete Bilevel Linear Programming}}.
\newblock {PhD}, Lehigh University, 2011.
\newblock URL \url{http://coral.ie.lehigh.edu/{~}ted/files/papers/
  ScottDeNegreDissertation11.pdf}.

\bibitem[DeNegre and Ralphs(2009)]{DeNRal09}
S.~DeNegre and T.K. Ralphs.
\newblock {A Branch-and-Cut Algorithm for Bilevel Integer Programming}.
\newblock In \emph{Proceedings of the Eleventh INFORMS Computing Society
  Meeting}, pages 65--78, 2009.
\newblock \doi{10.1007/978-0-387-88843-9_4}.
\newblock URL \url{http://coral.ie.lehigh.edu/~ted/files/papers/BILEVEL08.pdf}.

\bibitem[Dhamdhere et~al.(2005)Dhamdhere, Ravi, and Singh]{dhamdhere2005two}
K.~Dhamdhere, R.~Ravi, and M.~Singh.
\newblock On two-stage stochastic minimum spanning trees.
\newblock In \emph{International Conference on Integer Programming and
  Combinatorial Optimization}, pages 321--334. Springer, 2005.

\bibitem[Dyer and Stougie(2006)]{Dyer2006}
M.~Dyer and L.~Stougie.
\newblock Computational complexity of stochastic programming problems.
\newblock \emph{Mathematical Programming}, 106\penalty0 (3):\penalty0 423--432,
  May 2006.
\newblock ISSN 1436-4646.
\newblock \doi{10.1007/s10107-005-0597-0}.
\newblock URL \url{https://doi.org/10.1007/s10107-005-0597-0}.

\bibitem[Fa{\'\i}sca et~al.(2007)Fa{\'\i}sca, Dua, Rustem, Saraiva, and
  Pistikopoulos]{faisca-etal07}
N.P. Fa{\'\i}sca, V.~Dua, B.~Rustem, P.M. Saraiva, and E.N. Pistikopoulos.
\newblock Parametric global optimisation for bilevel programming.
\newblock \emph{Journal of Global Optimization}, 38:\penalty0 609--623, 2007.

\bibitem[Fischetti et~al.(2017)Fischetti, Ljubi{\'c}, Monaci, and
  Sinnl]{fischettietal17b}
M.~Fischetti, I.~Ljubi{\'c}, M.~Monaci, and M.~Sinnl.
\newblock A new general-purpose algorithm for mixed-integer bilevel linear
  programs.
\newblock \emph{Operations Research}, 65\penalty0 (6):\penalty0 1615--1637,
  2017.

\bibitem[Fischetti et~al.(2018)Fischetti, Ljubi{\'c}, Monaci, and
  Sinnl]{fischettietal17a}
M.~Fischetti, I.~Ljubi{\'c}, M.~Monaci, and M.~Sinnl.
\newblock On the use of intersection cuts for bilevel optimization.
\newblock \emph{Mathematical Programming}, 72:\penalty0 77--103, 2018.

\bibitem[Furini et~al.(2019)Furini, Ljubic, Martin, and Segundo]{furinimaximum}
F.~Furini, I.~Ljubic, S.~Martin, and P.~San Segundo.
\newblock The maximum clique interdiction game.
\newblock \emph{European Journal of Operational Research}, 277(1):\penalty0
  112--127, 2019.

\bibitem[Gade et~al.(2012)Gade, K{\"u}{\c{c}}{\"u}kyavuz, and
  Sen]{gade2012decomposition}
D.~Gade, S.~K{\"u}{\c{c}}{\"u}kyavuz, and S.~Sen.
\newblock Decomposition algorithms with parametric gomory cuts for two-stage
  stochastic integer programs.
\newblock \emph{Mathematical Programming}, pages 1--26, 2012.

\bibitem[Gan et~al.(2018)Gan, Elkind, and Wooldridge]{gan2018stackelberg}
J.~Gan, E.~Elkind, and M.~Wooldridge.
\newblock Stackelberg security games with multiple uncoordinated defenders.
\newblock In \emph{Proc. of 17th Int. Conf. on Autonomous Agents and Multiagent
  Systems (AAMAS 2008)}, 2018.

\bibitem[{Garc\'{e}s} et~al.(2009){Garc\'{e}s}, {Conejo},
  {Garc\'{i}a-Bertrand}, and {Romero}]{garcesetal09}
L.P. {Garc\'{e}s}, A.J. {Conejo}, R.~{Garc\'{i}a-Bertrand}, and R.~{Romero}.
\newblock A bilevel approach to transmission expansion planning within a market
  environment.
\newblock \emph{IEEE Transactions on Power Systems}, 24\penalty0 (3):\penalty0
  1513--1522, Aug 2009.

\bibitem[Garey and Johnson(1979)]{garey79}
M.R. Garey and D.S. Johnson.
\newblock \emph{Computers and Intractability: A Guide to the Thoery of
  NP-Completeness}.
\newblock W.H. Freeman and Company, 1979.

\bibitem[Gendreau et~al.(1996)Gendreau, Laporte, and
  S{\'e}guin]{gendreau1996stochastic}
M.~Gendreau, G.~Laporte, and R.~S{\'e}guin.
\newblock Stochastic vehicle routing.
\newblock \emph{European Journal of Operational Research}, 88\penalty0
  (1):\penalty0 3--12, 1996.

\bibitem[Ghare et~al.(1971)Ghare, Montgomery, and Turner]{ghare71}
P.M. Ghare, D.C. Montgomery, and W.C. Turner.
\newblock Optimal interdiction policy for a flow network.
\newblock \emph{Naval Research Logistics Quarterly}, 18:\penalty0 27--45, 1971.

\bibitem[G{\o}rtz et~al.(2012)G{\o}rtz, Nagarajan, and
  Saket]{gortz2012stochastic}
I.L. G{\o}rtz, V.~Nagarajan, and R.~Saket.
\newblock Stochastic vehicle routing with recourse.
\newblock In \emph{International Colloquium on Automata, Languages, and
  Programming}, pages 411--423. Springer, 2012.

\bibitem[Grass and Fischer(2016)]{grass2016two}
E.~Grass and K.~Fischer.
\newblock Two-stage stochastic programming in disaster management: A literature
  survey.
\newblock \emph{Surveys in Operations Research and Management Science},
  21\penalty0 (2):\penalty0 85--100, 2016.

\bibitem[Grimm et~al.(2016)Grimm, Martin, Martin, Weibelzahl, and
  Zöttl]{GRIMM2016493}
V.~Grimm, A.~Martin, M.~Martin, M.~Weibelzahl, and G.~Zöttl.
\newblock Transmission and generation investment in electricity markets: The
  effects of market splitting and network fee regimes.
\newblock \emph{European Journal of Operational Research}, 254\penalty0
  (2):\penalty0 493 -- 509, 2016.
\newblock ISSN 0377-2217.
\newblock \doi{https://doi.org/10.1016/j.ejor.2016.03.044}.
\newblock URL \url{http://www.sciencedirect.com/science/article/pii/
  S0377221716301904}.

\bibitem[Gupta et~al.(2007)Gupta, Ravi, and Sinha]{gupta2007lp}
A.~Gupta, R.~Ravi, and A.~Sinha.
\newblock Lp rounding approximation algorithms for stochastic network design.
\newblock \emph{Mathematics of Operations Research}, 32\penalty0 (2):\penalty0
  345--364, 2007.

\bibitem[G{\"{u}}zelsoy(2009)]{Guzelsoy2009}
M.~G{\"{u}}zelsoy.
\newblock \emph{{Dual Methods in Mixed Integer Linear Programming}}.
\newblock {PhD}, Lehigh University, 2009.
\newblock URL \url{http://coral.ie.lehigh.edu/{~}ted/files/papers/
  MenalGuzelsoyDissertation09.pdf}.

\bibitem[G{\"u}zelsoy and Ralphs(2007)]{GuzRal07}
M.~G{\"u}zelsoy and T.K. Ralphs.
\newblock {Duality for Mixed-Integer Linear Programs}.
\newblock \emph{International Journal of Operations Research}, 4:\penalty0
  118--137, 2007.
\newblock URL \url{http://coral.ie.lehigh.edu/~ted/files/papers/MILPD06.pdf}.

\bibitem[Hanasusanto et~al.(2016)Hanasusanto, Kuhn, and
  Wiesemann]{Hanasusanto2016}
G.A. Hanasusanto, D.~Kuhn, and W.~Wiesemann.
\newblock A comment on ``computational complexity of stochastic programming
  problems''.
\newblock \emph{Mathematical Programming}, 159\penalty0 (1):\penalty0 557--569,
  Sep 2016.
\newblock ISSN 1436-4646.
\newblock \doi{10.1007/s10107-015-0958-2}.
\newblock URL \url{https://doi.org/10.1007/s10107-015-0958-2}.

\bibitem[Hansen et~al.(1992)Hansen, Jaumard, and Savard]{hansen92}
P.~Hansen, B.~Jaumard, and G.~Savard.
\newblock New branch-and-bound rules for linear bilevel programming.
\newblock \emph{SIAM Journal on Scientific and Statistical Computing},
  13\penalty0 (5):\penalty0 1194--1217, 1992.

\bibitem[Hassanzadeh and Ralphs(2014{\natexlab{a}})]{HasRal14}
A.~Hassanzadeh and T.K. Ralphs.
\newblock {A Generalized Benders' Algorithm for Two-Stage Stochastic Program
  with Mixed Integer Recourse}.
\newblock Technical Report COR@L Laboratory Report 14T-005, Lehigh University,
  2014{\natexlab{a}}.
\newblock URL \url{http://coral.ie.lehigh.edu/~ted/files/papers/
  SMILPGenBenders14.pdf}.

\bibitem[Hassanzadeh and Ralphs(2014{\natexlab{b}})]{HasRal14-1}
A.~Hassanzadeh and T.K. Ralphs.
\newblock {On the Value Function of a Mixed Integer Linear Optimization Problem
  and an Algorithm for Its Construction}.
\newblock Technical report, COR@L Laboratory Report 14T-004, Lehigh University,
  2014{\natexlab{b}}.
\newblock URL \url{http://coral.ie.lehigh.edu/~ted/files/papers/
  MILPValueFunction14.pdf}.

\bibitem[Held and Woodruff(2005)]{held05}
H.~Held and D.L. Woodruff.
\newblock Heuristics for multi-stage interdiction of stochastic networks.
\newblock \emph{Journal of Heuristics}, 11\penalty0 (5-6):\penalty0 483--500,
  2005.

\bibitem[Hemmati and Smith(2016)]{HemSmi16}
M.~Hemmati and J.C. Smith.
\newblock A mixed integer bilevel programming approach for a competitive set
  covering problem.
\newblock Technical report, Clemson University, 2016.

\bibitem[Hobbs and Nelson(1992)]{hobbs-nelson92}
B.F. Hobbs and S.K. Nelson.
\newblock A nonlinear bilevel model for analysis of electric utility
  demand-side planning issues.
\newblock \emph{Annals of Operations Research}, 34\penalty0 (1):\penalty0
  255--274, 1992.

\bibitem[Israeli(1999)]{israeli99}
E.~Israeli.
\newblock \emph{System Interdiction and Defense}.
\newblock PhD thesis, Naval Postgraduate School, 1999.

\bibitem[Israeli and Wood(2002)]{israeli02}
E.~Israeli and R.K. Wood.
\newblock Shortest path network interdiction.
\newblock \emph{Networks}, 40\penalty0 (2):\penalty0 97--111, 2002.

\bibitem[Janjarassuk and Linderoth(2008)]{janjarassuk06}
U.~Janjarassuk and J.~Linderoth.
\newblock Reformulation and sampling to solve a stochastic network interdiction
  problem.
\newblock \emph{Networks}, 52:\penalty0 120--132, 2008.

\bibitem[Jeroslow(1985)]{Jeroslow1985}
R.G. Jeroslow.
\newblock The polynomial hierarchy and a simple model for competitive analysis.
\newblock \emph{Mathematical Programming}, 32\penalty0 (2):\penalty0 146--164,
  Jun 1985.
\newblock ISSN 1436-4646.
\newblock \doi{10.1007/BF01586088}.
\newblock URL \url{https://doi.org/10.1007/BF01586088}.

\bibitem[Kall and Mayer(2010)]{kall2010stochastic}
P.~Kall and J.~Mayer.
\newblock \emph{{Stochastic linear programming: models, theory, and
  computation}}.
\newblock Springer Verlag, 2010.
\newblock ISBN 1441977287.

\bibitem[Kara and Verter(2004)]{kara-verter04}
B.~Kara and V.~Verter.
\newblock Designing a road network for hazardous materials transportation.
\newblock \emph{Transportation Science}, 38\penalty0 (2):\penalty0 188--196,
  2004.

\bibitem[Karp(1975)]{karp1975computational}
R.M. Karp.
\newblock On the computational complexity of combinatorial problems.
\newblock \emph{Networks}, 5:\penalty0 45--68, 1975.

\bibitem[Katriel et~al.(2007)Katriel, Kenyon-Mathieu, and
  Upfal]{katriel2007commitment}
I.~Katriel, C.~Kenyon-Mathieu, and E.~Upfal.
\newblock Commitment under uncertainty: Two-stage stochastic matching problems.
\newblock In \emph{International Colloquium on Automata, Languages, and
  Programming}, pages 171--182. Springer, 2007.

\bibitem[Kiekintveld et~al.(2009)Kiekintveld, Jain, Tsai, Pita,
  Ord{\'{o}}{\~{n}}ez, and Tambe]{KiekintveldJTPOT09}
C.~Kiekintveld, M.~Jain, J.~Tsai, J.~Pita, F.~Ord{\'{o}}{\~{n}}ez, and
  M.~Tambe.
\newblock Computing optimal randomized resource allocations for massive
  security games.
\newblock In \emph{AAMAS}, pages 689--696, 2009.

\bibitem[Klein(2019)]{klein2019complexity}
K.-M. Klein.
\newblock About the complexity of two-stage stochastic ips, 2019.

\bibitem[Kong et~al.(2006)Kong, Schaefer, and Hunsaker]{kong2006two}
N.~Kong, A.J. Schaefer, and B.~Hunsaker.
\newblock {Two-stage integer programs with stochastic right-hand sides: a
  superadditive dual approach}.
\newblock \emph{Mathematical Programming}, 108\penalty0 (2):\penalty0 275--296,
  2006.

\bibitem[K{\"{o}}ppe et~al.(2010)K{\"{o}}ppe, Queyranne, and Ryan]{koppe10}
M.~K{\"{o}}ppe, M.~Queyranne, and C.~T. Ryan.
\newblock Parametric integer programming algorithm for bilevel mixed integer
  programs.
\newblock \emph{Journal of Optimization Theory and Applications}, 146\penalty0
  (1):\penalty0 137--150, Jul 2010.

\bibitem[Kulkarni and Shanbhag(2014)]{kulkarni2014shared}
A.A. Kulkarni and U.V. Shanbhag.
\newblock A shared-constraint approach to multi-leader multi-follower games.
\newblock \emph{Set-Valued and Variational Analysis}, 22\penalty0 (4):\penalty0
  691--720, 2014.

\bibitem[Labb{\'e} and Violin(2013)]{labbe2013bilevel}
M.~Labb{\'e} and A.~Violin.
\newblock Bilevel programming and price setting problems.
\newblock \emph{4OR}, 11\penalty0 (1):\penalty0 1--30, 2013.

\bibitem[Labb\'{e} et~al.(1998)Labb\'{e}, Marcotte, and Savard]{labbe98}
M.~Labb\'{e}, P.~Marcotte, and G.~Savard.
\newblock A bilevel model of taxation and its application to optimal highway
  pricing.
\newblock \emph{Management Science}, 44:\penalty0 1608--1622, 1998.

\bibitem[Laporte and Louveaux(1993)]{laporte1993integer}
G.~Laporte and F.V. Louveaux.
\newblock The integer l-shaped method for stochastic integer programs with
  complete recourse.
\newblock \emph{Operations research letters}, 13\penalty0 (3):\penalty0
  133--142, 1993.

\bibitem[Laszka et~al.(2016)Laszka, Lou, and Vorobeychik]{laszka2016multi}
A.~Laszka, J.~Lou, and Y.~Vorobeychik.
\newblock Multi-defender strategic filtering against spear-phishing attacks.
\newblock In \emph{Proc. of 30th AAAI Conf. on Artificial Intelligence (AAAI
  2016)}, 2016.

\bibitem[Leyffer and Munson(2010)]{leyffer2010solving}
S.~Leyffer and T.~Munson.
\newblock Solving multi-leader-common-follower games.
\newblock \emph{OPT MET SO}, 25\penalty0 (4):\penalty0 601--623, 2010.

\bibitem[Lodi et~al.(2009)Lodi, Ralphs, Rossi, and
  Smriglio]{lodi2009interdiction}
A.~Lodi, T.K. Ralphs, F.~Rossi, and S.~Smriglio.
\newblock Interdiction branching.
\newblock Technical Report OR/09/10, DEIS-Universit{\`a} di Bologna, 2009.

\bibitem[Lou and Vorobeychik(2015)]{lou2015equilibrium}
J.~Lou and Y.~Vorobeychik.
\newblock Equilibrium analysis of multi-defender security games.
\newblock In \emph{Proc. of 24th Int. Joint Conf. on Artificial Intelligence
  (IJCAI 2019)}, 2015.

\bibitem[Lou et~al.(2017)Lou, Smith, and Vorobeychik]{lou2017multidefender}
J.~Lou, A.M. Smith, and Y.~Vorobeychik.
\newblock Multidefender security games.
\newblock \emph{IEEE INTELL SYST}, 32\penalty0 (1):\penalty0 50--60, 2017.

\bibitem[Louveaux and van~der Vlerk(1993)]{louveaux1993stochastic}
F.V. Louveaux and M.H. van~der Vlerk.
\newblock {Stochastic programming with simple integer recourse}.
\newblock \emph{Mathematical programming}, 61\penalty0 (1):\penalty0 301--325,
  1993.

\bibitem[Lozano and Smith(2017)]{lozanosmith17}
L.~Lozano and J.C. Smith.
\newblock A value-function-based exact approach for the bilevel mixed-integer
  programming problem.
\newblock \emph{Operations Research}, 65\penalty0 (3):\penalty0 768--786, 2017.

\bibitem[Mahajan(2009)]{Mahajan2009}
A.~Mahajan.
\newblock \emph{{On Selecting Disjunctions in Mixed Integer Linear
  Programming}}.
\newblock {PhD}, Lehigh University, 2009.
\newblock URL \url{http://coral.ie.lehigh.edu/{~}ted/files/papers/
  AshutoshMahajanDissertation09.pdf}.

\bibitem[Mahajan and Ralphs(2010)]{MahRal10}
A.~Mahajan and T.K. Ralphs.
\newblock {On the Complexity of Selecting Disjunctions in Integer Programming}.
\newblock \emph{SIAM Journal on Optimization}, 20\penalty0 (5):\penalty0
  2181--2198, 2010.
\newblock \doi{10.1137/080737587}.
\newblock URL \url{http://coral.ie.lehigh.edu/~ted/files/papers/
  Branching08.pdf}.

\bibitem[Marchesi et~al.(2018)Marchesi, Coniglio, and
  Gatti]{marchesi2018leadership}
A.~Marchesi, S.~Coniglio, and N.~Gatti.
\newblock Leadership in singleton congestion games.
\newblock In \emph{Proc. of 27th Int. Joint Conf. on Artificial Intelligence
  (IJCAI 2018)}, pages 447--453. AAAI Press, 2018.

\bibitem[McMasters and Mustin(1970)]{mcmasters70}
A.W. McMasters and T.M. Mustin.
\newblock Optimal interdiction of a supply network.
\newblock \emph{Naval Research Logistics Quarterly}, 17:\penalty0 261--268,
  1970.

\bibitem[Migdalas(1995)]{migdalas95}
A.~Migdalas.
\newblock Bilevel programming in traffic planning: models, methods and
  challenge.
\newblock \emph{Journal of Global Optimization}, 7:\penalty0 381--405, 1995.

\bibitem[Moore and Bard(1990)]{moore90}
J.T. Moore and J.F. Bard.
\newblock The mixed integer linear bilevel programming problem.
\newblock \emph{Operations Research}, 38\penalty0 (5):\penalty0 911--921, 1990.

\bibitem[Pang and Fukushima(2005)]{pang2005quasi}
J-S. Pang and M.~Fukushima.
\newblock Quasi-variational inequalities, generalized nash equilibria, and
  multi-leader-follower games.
\newblock \emph{Computational Management Science}, 2\penalty0 (1):\penalty0
  21--56, 2005.

\bibitem[Paruchuri et~al.(2008)Paruchuri, Pearce, Marecki, Tambe, Ordonez, and
  Kraus]{paruchuri2008playing}
P.~Paruchuri, J.P. Pearce, J.~Marecki, M.~Tambe, F.~Ordonez, and S.~Kraus.
\newblock Playing games for security: an efficient exact algorithm for solving
  bayesian stackelberg games.
\newblock In \emph{Proc. of 7th Int. Conf. on Autonomous Agents and Multiagent
  Systems (AAMAS 2008)}, pages 895--902, 2008.

\bibitem[Ralphs and G{\"u}zelsoy(2005)]{RalGuz05}
T.K. Ralphs and M.~G{\"u}zelsoy.
\newblock {The SYMPHONY Callable Library for Mixed Integer Programming}.
\newblock In \emph{Proceedings of the Ninth INFORMS Computing Society
  Conference}, pages 61--76, 2005.
\newblock \doi{10.1007/0-387-23529-9_5}.
\newblock URL \url{http://coral.ie.lehigh.edu/~ted/files/papers/
  SYMPHONY04.pdf}.

\bibitem[Ralphs and G{\"u}zelsoy(2006)]{RalGuz06}
T.K. Ralphs and M.~G{\"u}zelsoy.
\newblock {Duality and Warm Starting in Integer Programming}.
\newblock In \emph{The Proceedings of the 2006 NSF Design, Service, and
  Manufacturing Grantees and Research Conference}, 2006.
\newblock URL \url{http://coral.ie.lehigh.edu/~ted/files/papers/DMII06.pdf}.

\bibitem[Rockafellar and Uryasev(2000)]{rockafellar2000optimization}
R.T. Rockafellar and S.~Uryasev.
\newblock Optimization of conditional value-at-risk.
\newblock \emph{{Journal of Risk}}, 2:\penalty0 21--42, 2000.

\bibitem[Rutenburg(1994)]{rutenburg1994propositional}
V.~Rutenburg.
\newblock Propositional truth maintenance systems: Classification and
  complexity analysis.
\newblock \emph{Annals of Mathematics and Artificial Intelligence}, 10\penalty0
  (3):\penalty0 207--231, 1994.

\bibitem[Saharidis and Ierapetritou(2008)]{saharidis-ierapetritou08}
G.K. Saharidis and M.G. Ierapetritou.
\newblock Resolution method for mixed integer bi-level linear problems based on
  decomposition technique.
\newblock \emph{Journal of Global Optimization}, 44\penalty0 (1):\penalty0
  29--51, 2008.

\bibitem[Sandholm(2002)]{sandholm2002evolutionary}
W.H. Sandholm.
\newblock Evolutionary implementation and congestion pricing.
\newblock \emph{The Review of Economic Studies}, 69\penalty0 (3):\penalty0
  667--689, 2002.

\bibitem[Scaparra and Church(2008)]{ScaparraChurch08COR}
M.P. Scaparra and R.L. Church.
\newblock A bilevel mixed-integer program for critical infrastructure
  protection planning.
\newblock \emph{Computers and Operations Research}, 35\penalty0 (6):\penalty0
  1905--1923, 2008.
\newblock ISSN 0305-0548.

\bibitem[Schaefer and Umans(2002)]{schaefer2002completeness}
M.~Schaefer and C.~Umans.
\newblock Completeness in the polynomial-time hierarchy: A compendium.
\newblock \emph{SIGACT news}, 33\penalty0 (3):\penalty0 32--49, 2002.

\bibitem[Schultz et~al.(1998)Schultz, Stougie, and Van
  Der~Vlerk]{schultz1998solving}
R.~Schultz, L.~Stougie, and M.H. Van Der~Vlerk.
\newblock {Solving stochastic programs with integer recourse by enumeration: A
  framework using Gr{\"o}bner basis}.
\newblock \emph{Mathematical Programming}, 83\penalty0 (1):\penalty0 229--252,
  1998.

\bibitem[Segundo et~al.(2019)Segundo, Coniglio, Furini, and
  Ljubi{\'c}]{san2019new}
P.~San Segundo, S.~Coniglio, F.~Furini, and I.~Ljubi{\'c}.
\newblock A new branch-and-bound algorithm for the maximum edge-weighted clique
  problem.
\newblock \emph{European Journal of Operational Research}, 278\penalty0
  (1):\penalty0 76--90, 2019.

\bibitem[Sen and Higle(2005)]{sen2005c}
S.~Sen and J.L. Higle.
\newblock {The {$C^3$} theorem and a {$D^2$} algorithm for large scale
  stochastic mixed-integer programming: Set convexification}.
\newblock \emph{Mathematical Programming}, 104\penalty0 (1):\penalty0 1--20,
  2005.
\newblock ISSN 0025-5610.

\bibitem[Shapiro(2003)]{shapiro2003monte}
A.~Shapiro.
\newblock Monte carlo sampling methods.
\newblock \emph{Handbooks in operations research and management science},
  10:\penalty0 353--425, 2003.

\bibitem[Sherali and Fraticelli(2002)]{sherali2002modification}
H.D. Sherali and B.M.P. Fraticelli.
\newblock {A modification of {B}enders' decomposition algorithm for discrete
  subproblems: An approach for stochastic programs with integer recourse}.
\newblock \emph{Journal of Global Optimization}, 22\penalty0 (1):\penalty0
  319--342, 2002.

\bibitem[Sherali and Zhu(2006)]{sherali2006solving}
H.D. Sherali and X.~Zhu.
\newblock {On solving discrete two-stage stochastic programs having
  mixed-integer first-and second-stage variables}.
\newblock \emph{Mathematical Programming}, 108\penalty0 (2):\penalty0 597--616,
  2006.

\bibitem[Smith et~al.(2014)Smith, Vorobeychik, and
  Letchford]{smith2014multidefender}
A.~Smith, Y.~Vorobeychik, and J.~Letchford.
\newblock Multidefender security games on networks.
\newblock \emph{ACM SIGMETRICS Performance Evaluation Review}, 41\penalty0
  (4):\penalty0 4--7, 2014.

\bibitem[Stackelberg(2010)]{von2010market}
H.~Von Stackelberg.
\newblock \emph{Market structure and equilibrium}.
\newblock Springer Science \& Business Media, 2010.

\bibitem[Stockmeyer(1976)]{stockmeyer77}
L.J. Stockmeyer.
\newblock The polynomial-time hierarchy.
\newblock \emph{Theoretical Computer Science}, 3:\penalty0 1--22, 1976.

\bibitem[Tahernejad(2019)]{Tahernejad2019}
S.~Tahernejad.
\newblock \emph{{Two-stage Mixed Integer Stochastic Bilevel Linear
  Optimization}}.
\newblock {PhD}, Lehigh University, 2019.

\bibitem[Tahernejad et~al.(2016)Tahernejad, Ralphs, and DeNegre]{TahRalDeN16}
S.~Tahernejad, T.K. Ralphs, and S.T. DeNegre.
\newblock {A Branch-and-Cut Algorithm for Mixed Integer Bilevel Linear
  Optimization Problems and Its Implementation}.
\newblock \emph{Mathematical Programming Computation (to appear)}, 2016.
\newblock URL \url{http://coral.ie.lehigh.edu/~ted/files/papers/MIBLP16.pdf}.
\newblock To appear, Mathematical Programming Computation.

\bibitem[Tamble(2011)]{tambe2011security}
M.~Tamble.
\newblock \emph{Security and Game Theory: Algorithms, Deployed Systems, Lessons
  Learned}.
\newblock Cambridge University Press, 2011.
\newblock ISBN 1107096421, 9781107096424.

\bibitem[Uryasev(2000)]{uryasev2000conditional}
S.~Uryasev.
\newblock Conditional value-at-risk: Optimization algorithms and applications.
\newblock In \emph{Proceedings of the IEEE/IAFE/INFORMS 2000 Conference on
  Computational Intelligence for Financial Engineering (CIFEr)(Cat. No.
  00TH8520)}, pages 49--57. IEEE, 2000.

\bibitem[Van~Slyke and Wets(1969)]{van1969shaped}
R.M. Van~Slyke and R.~Wets.
\newblock L-shaped linear programs with applications to optimal control and
  stochastic programming.
\newblock \emph{SIAM Journal on Applied Mathematics}, pages 638--663, 1969.

\bibitem[Vicente et~al.(1996)Vicente, Savard, and J{\'u}dice]{vicente96}
L.~Vicente, G.~Savard, and J.~J{\'u}dice.
\newblock Discrete linear bilevel programming problem.
\newblock \emph{Journal of Optimization Theory and Applications}, 89\penalty0
  (3):\penalty0 597--614, 1996.

\bibitem[Vicente and Calamai(1994)]{vicente-calamai94}
L.N. Vicente and P.H. Calamai.
\newblock Bilevel and multilevel programming: A bibliography review.
\newblock \emph{Journal of Global Optimization}, 5\penalty0 (3):\penalty0
  291--306, 1994.

\bibitem[von Stengel and Zamir(2010)]{von2010leadership}
B.~von Stengel and S.~Zamir.
\newblock Leadership games with convex strategy sets.
\newblock \emph{Games and Economic Behavior}, 69\penalty0 (2):\penalty0
  446--457, 2010.

\bibitem[Wallace and Ziemba(2005)]{wallace2005applications}
S.W. Wallace and W.T. Ziemba.
\newblock \emph{Applications of stochastic programming}.
\newblock SIAM, 2005.

\bibitem[Wang and Xu(2017)]{wangxu17}
L.~Wang and P.~Xu.
\newblock The watermelon algorithm for the bilevel integer linear programming
  problem.
\newblock \emph{SIAM Journal on Optimization}, 27\penalty0 (3):\penalty0
  1403--1430, 2017.

\bibitem[Wen and Yang(1990)]{wen90}
U.P. Wen and Y.H. Yang.
\newblock Algorithms for solving the mixed integer two-level linear programming
  problem.
\newblock \emph{Computers \& Operations Research}, 17\penalty0 (2):\penalty0
  133--142, 1990.

\bibitem[Williams(1996)]{williams96-2}
H.P. Williams.
\newblock Duality in mathematics and linear and integer programming.
\newblock \emph{Journal of Optimization Theory and Applications}, 90\penalty0
  (2):\penalty0 257--278, 1996.

\bibitem[Wollmer(1964)]{wollmer64}
R.~Wollmer.
\newblock Removing arcs from a network.
\newblock \emph{Operations Research}, 12\penalty0 (6):\penalty0 934--940, 1964.

\bibitem[Wolsey(1981{\natexlab{a}})]{wolsey1981integer}
L.A. Wolsey.
\newblock {Integer programming duality: Price functions and sensitivity
  analysis}.
\newblock \emph{Mathematical Programming}, 20\penalty0 (1):\penalty0 173--195,
  1981{\natexlab{a}}.
\newblock ISSN 0025-5610.

\bibitem[Wolsey(1981{\natexlab{b}})]{wolsey81}
L.A. Wolsey.
\newblock Integer programming duality: Price functions and sensitivity
  analaysis.
\newblock \emph{Mathematical Programming}, 20:\penalty0 173--195,
  1981{\natexlab{b}}.

\bibitem[Wood(1993)]{wood93}
R.K. Wood.
\newblock Deterministic network interdiction.
\newblock \emph{Mathematical and Computer Modelling}, 17\penalty0 (2):\penalty0
  1--18, 1993.

\bibitem[Xu and Wang(2014)]{xuwang14}
P.~Xu and L.~Wang.
\newblock An exact algorithm for the bilevel mixed integer linear programming
  problem under three simplifying assumptions.
\newblock \emph{Computers \& operations research}, 41:\penalty0 309--318, 2014.

\bibitem[Zeng and An(2014)]{zengan14}
B.~Zeng and U.~An.
\newblock Solving bilevel mixed integer program by reformulations and
  decomposition.
\newblock \emph{Optimization online}, pages 1--34, 2014.

\bibitem[Zhang et~al.(2016)Zhang, Snyder, Ralphs, and Xue]{ZhaSnyRalXue16}
Y.~Zhang, L.V. Snyder, T.K. Ralphs, and Z.~Xue.
\newblock {The Competitive Facility Location Problem Under Disruption Risks}.
\newblock \emph{Transportation Research Part E: Logistics and Transportation
  Review}, 93, 2016.
\newblock ISSN 13665545.
\newblock \doi{10.1016/j.tre.2016.07.002}.
\newblock URL \url{http://coral.ie.lehigh.edu/~ted/files/papers/CFLPD16.pdf}.

\end{thebibliography}
\end{document}